\documentclass[10pt, a4paper]{amsart}
\input xy   
\xyoption{all}
\numberwithin{equation}{section}

\swapnumbers

\theoremstyle{plain}

  \begingroup
        \newtheorem{theorem}[equation]{Theorem}
        
        \newtheorem{proposition}[equation]{Proposition}
        \newtheorem{corollary}[equation]{Corollary}
        \newtheorem{fact}[equation]{Fact}
        \newtheorem{facts}[equation]{Facts}
        \newtheorem{assumption}[equation]{Assumption}
        \newtheorem{remark}[equation]{Remark}
        \newtheorem{definition}[equation]{Definition}

        \newtheorem{notation}[equation]{Notation}

        \newtheorem{sinnada}[equation]{}
  \endgroup
\theoremstyle{definition}
  \begingroup
        \newtheorem{example}[equation]{Example}
        \newtheorem{yoneda}[equation]{Yoneda}
        
  \endgroup

\newcommand{\mr}[1]{\buildrel {#1} \over \longrightarrow}

\newcommand{\siff}{\Leftrightarrow}

\newcommand{\rimply}{\Rightarrow}

\newcommand{\mono}{\hookrightarrow}
\newcommand{\mmr}[1]{\buildrel {#1} \over \hookrightarrow}

\newcommand{\cc}{\mathcal}
\newcommand{\bb}{\mathbb}

\begin{document}

\title{Quasitopoi over a base category}

\author{Eduardo J. Dubuc  and  Luis Espa\~nol}

\begin{abstract}
In this paper we develop the theory of quasispaces (for a Grothendieck
topology) and of concrete quasitopoi, over a suitable base category. We introduce the notion of \emph{$f$-regular category} and of \emph{$f$-regular
functor}. The $f$-regular categories are
regular categories in which every family with a common codomain can be factorized into a
strict epimorphic family followed by a (single) monomorphism. The $f$-regular functors are (essentially)
functors that preserve finite strict monomorphic and arbitrary strict
epimorphic families.
These two concepts furnish the context to develop the 
constructions of the theory of concrete quasitopoi over a suitable
base category, which is a theory of \emph{pointed quasitopoi}. Our results on quasispaces and quasitopoi , or closely related ones,
were already established by Penon in \cite{P2}, but we prove them
here with different assumptions, and under a completely different light.  
\end{abstract}

\maketitle

{\sc introduction} \indent \vspace{1ex}

In this paper we develop the theory of quasispaces (for a Grothendieck
topology) and of
concrete quasitopoi, over a suitable base category. We use the
systematic theory of families of arrows in a category 
(including a characterization of strict epimorphic families as final
surjective families in a general context) \mbox{that we set in
  \cite{DE}.}

We introduce the notion of \emph{$f$-regular category} and of \emph{$f$-regular
functor}. The $f$-regular categories are
regular categories which have $f$-factorizations, a  non elementary
cocompletness condition:
\emph{A category has $f$-factorizations if every family
 $g_\alpha:X_\alpha\to X$ factorizes in the form $g_\alpha=m\circ
 h_\alpha$,
 where $m:H\to X$ is a monomorphism and $h_\alpha:X_\alpha\to H$
is an 
strict epimorphic family.} 
The $f$-regular functors are (essentially)
functors that preserve finite strict monomorphic and arbitrary strict
epimorphic families.
These two concepts furnish the context to develop the 
constructions of the theory of concrete quasitopoi.

Our results on quasispaces and quasitopoi , or closely related ones,
were already established by Penon in \cite{P2}, but we prove them
here with different assumptions, and under a completely different light.  
The notion of bounded quasitopos \cite{P2} leaves out the 
paradigmatic examples of quasitopoi, namely, 
the \emph{legitimate} categories of separated sheaves over a  
\emph{large} site. In this case, the categories
of sheaves are not topoi, and they are not even legitimate categories
since the class of morphisms between two sheaves is not in general a
set. However, separate sheaves form legitimate categories 
(insofar the
hom-sets are small) which are elementary quasitopoi (in the sense of
Penon). They  bear to elementary quasitopoi a relation which should 
be considered as corresponding (in the context of quasitopoi) to the
relation that Grothendieck topoi bear to elementary 
topoi. In this context, the size condition (\emph{bounded}) is
unnecessary, and probably misleading. We introduce the abstract notion of
$f$-quasitopos, which describes this situation, and generalize the
concrete quasitopoi of \cite{D}. Our notion is intermediate: 
$$ 
\text{bounded quasitopos} \;\;  \rimply  \;\; f\text{-quasitopos} \;\;\rimply 
\;\;\text{elementary quasitopos}
$$ 
The morphisms between $f$-quasitopoi are the $f$-regular functors,
which correspond to inverse images of geometric morphisms. The category
of $f$-quasitopoi over a suitable base category is a theory of
\emph{pointed} quasitopoi. 

\tableofcontents

\section{Families of arrows and topological functors}  \label{families}

\vspace{1ex}

In this section we recall briefly some notions and results from \cite{DE} that
we shall explicitly need, and in this way fix notation and
terminology. For comments and proofs we refer to \cite{DE}.

\vspace{1ex}

Given a category $\mathbb{T}$ and an object $X$ in $\mathbb{T}$, we 
shall work with families 
\mbox{$(X_{\alpha} \mr{g_{\alpha}} X)_{\alpha \in \Gamma}$} of 
arrows of $\mathbb{T}$ with codomain $X$. 

\begin{notation}
Given a family $(X_{\alpha} \mr{g_{\alpha}} X)_{\alpha \in \Gamma}$,
    we shall simply write $X_{\alpha} \mr{g_{\alpha}} X$, 
    omitting as well a label for the index set (the context will
    always tell whether we are considering a single $\alpha$ or the
    whole family). 

The diagrammatic notation always denotes a commutative diagram, unless
otherwise explicitly indicated.
\end{notation}
It is important to point 
out that we allow the families to be \emph{large}, that is, not
indexed by a set.

\vspace{1ex}

\begin{definition} \label{refinement}
We say that a family $Y_\lambda  \to X$ \emph{refines (is a refinement
 of)}  a
family $X_\alpha \to  X$ if there is a function between the indices 
$\lambda \mapsto \alpha_\lambda$ together with arrows 
$Y_\lambda \to X_{\alpha_\lambda}$ such that 
$$
\xymatrix@R=15pt@C=7pt
                 {
                  Y_\lambda  \ar[rr]  \ar[rd] &&  X_{\alpha_\lambda} \ar[ld]
                  \\
                  & X
                 }
$$
\end{definition}

\begin{definition} \label{rpullback}
Given an arrow $Y \to X$, we say than a family $Y_\lambda  \to Y$ is a 
\mbox{$r$\emph{-pull-back}} of a family  $X_\alpha \to  X$, if 
$Y_\lambda  \to Y \to X$ refines $X_\alpha \to  X$. That is, if there is a function between the indices $\lambda \mapsto \alpha_\lambda$ together with arrows 
$Y_\lambda \to X_{\alpha_\lambda}$ such that 
$$
\xymatrix@R=15pt@C=20pt
        {
         Y_{\lambda} \ar[r] \ar[d] 
        & X_{\alpha_{\lambda}} \ar[d]
        \\
        Y \ar[r] 
        & X
        }
$$
\end{definition}

We consider collections $\cc{A}$ of classes of 
families of arrows with common codomain, one class $\cc{A}_X$ 
(eventually empty) for each object $X$ in $\mathbb{T}$. We say that a
 family in $\cc{A}_X$ is a $\cc{A}$-family over $X$. 

\vspace{1ex}

\begin{definition}[operations on collections] \label{familyoperations} 
$ $ 

(1) We denote by $\cc{I}so$ the collection whose only arrows are the
isomorphisms.

\vspace{1ex}

(2) Given two collections $\cc{A}$, $\cc{B}$ we  define the \emph{composite} 
$ \cc{C} = \cc{A} \circ \cc{B}$ by means of the following implication:  
$$
  X_{\alpha} \rightarrow X \; \in \; \cc{A}_X \;\; and \;\;
  \forall \, \alpha \;\; X_{\alpha,\, \beta} \rightarrow X_{\alpha}
\; \in \; 
  \cc{B}_{X_\alpha}  
  \Longrightarrow \;\;
  X_{\alpha, \, \beta} \rightarrow X_{\alpha} \rightarrow X \;
  \in \; \cc{C}_X
$$
 
\vspace*{4ex}

(3) Given $\cc{A}$ we define a new collection, denoted $\pi \cc{A}$, by: 
\vspace{1ex}

$Y_{\alpha} \rightarrow Y \in \pi \cc{A}\;\; \iff \;\;$ there is
$X_\alpha \to X  \in \cc{A}$ and $\;Y\to X$ such that:
$$
\xymatrix@R=5pt{\\ the \; squares\;\;}
\xymatrix@R=15pt@C=20pt
              {
                Y_\alpha \ar[r] \ar[d] 
                & X_\alpha \ar[d]
                \\
                Y \ar[r] 
               & X
              }
\xymatrix@R=5pt{\\ \;\; are \; pullbacks \; for \; all \; \alpha.}
$$

(4) Given $\cc{A}$ we define a new collection, denoted $s \cc{A}$, by:

\vspace{1ex}

$X_\alpha  \to X \;\in\; s \cc{A} \;\; \iff \;\;$ there is a
refinement  
by a family $Y_\lambda \to X \; \in \cc{A}_X$.
\end{definition}                        

Notice that $\cc{A}\,\subseteq\, \pi \cc{A}$, and $\cc{A}\,\subseteq\, s \cc{A}$. We set now some properties of collections $\cc{A}$ defined by means of these operations:

\begin{definition}[properties of collections] \label{familyproperties}
$ $   
 \vspace{1ex}
 
(I) \emph{Isomorphisms}:  $\cc{I}so \,\subseteq \, \cc{A}.$ 

\vspace{1ex}

(C) \emph{Closed under composition}: 
                            $\cc{A} \circ \cc{A} \,\subseteq\, \cc{A}.$

\vspace{1ex}

(S) \emph{Saturated}: $s \cc{A} \,\subseteq\, \cc{A}$ (hence $\cc{A}=s \cc{A}$).

\vspace{1ex}

(U) \emph{Universal}: Given $X_\alpha \to X \, \in \cc{A}$ and $Y \to X$, there exists an r-pull-back $Y_\lambda \to Y \, \in \cc{A}$:
$$
\xymatrix@R=15pt@C=20pt
        {
         Y_{\lambda} \ar[r] \ar[d] 
        & X_{\alpha_{\lambda}} \ar[d]
        \\
        Y \ar[r] 
        & X
        }
$$

\vspace{1ex}

If  a collection $\cc{A}$ satisfies (S), and the category has finite limits, then (U) is 
equivalent to:

\vspace{1ex}

(U) \emph{Stable under pullback}: $\pi\cc{A} \subseteq \cc{A}$ (hence
$\cc{A}=\pi \cc{A}$).

\vspace{1ex}

(F) \emph{Filtered}: Given 
$X_{\alpha}\to X \, \in \, \cc{A}_X,\;  Y_{\beta}\to X \, \in \,
  \cc{A}_X$, there exists a common refinement 
$Z_{\lambda}\to X \; \in \; \cc{A}_X$:  
$$
\xymatrix@R=15pt@C=20pt
         {
             Z_{\lambda}  \ar[d]  \ar[r] \ar[rd]
            &  X_{\alpha_{\lambda}} \ar[d] 
          \\
          Y_{\beta_{\lambda}} \ar[r]  &  X 
         }
$$ 
\end{definition}

\begin{facts} [about collections]  \label{collectionfacts} 
The following statements about given collections hold:

\begin{enumerate} 

\item  \label{UCimpliesF}
If a collection satisfies the properties (U) and (C) then it satisfies (F).

\item If $\cc{A}$ and $\cc{B}$ both satisfy (I), (resp. (C)),
  (resp. (S)), then so does $\cc{A} \cap \cc{B}$.

\item  \label{intersection}
If  $\cc{A}$ and $\cc{B}$ both satisfy (S) and (U), (resp. (S)
  and (F)), then so does  $\cc{A} \cap \cc{B}$.
\end{enumerate}
\end{facts} 

An important collection of families are the strict epimorphic
families. We recall now this notion from SGA4 \cite[I, 10.3, p. 180]{G2}:  

\begin{definition} \label{strictepi}
Given two families of arrows 
$f_{\alpha}:X_{\alpha}\to X$,  $g_{\alpha}:X_{\alpha}\to Y$, 
with the same indexes and domains, we say that  $g_{\alpha}$ is 
\emph{compatible with} $f_{\alpha}$ if for any pair of arrows 
$(x_{\alpha}:Z\to X_{\alpha}, \; x_{\beta}:Z\to X_{\beta})$ with the 
same domain the following condition holds: 
$f_{\alpha}\circ x_{\alpha}= f_{\beta}\circ x_{\beta}$ implies 
$g_{\alpha}\circ x_{\alpha}= g_{\beta}\circ x_{\beta}$.
 
A family $f_{\alpha}:X_{\alpha}\to X$ is \emph{strict epimorphic} 
if for any family $g_{\alpha}:X_{\alpha}\to Y$ \mbox{which is} 
compatible with $f_{\alpha}$, there exists a unique $g:X\to Y$ 
such that $g\circ f_{\alpha}=g_{\alpha}$ for all $\alpha$. 
\end{definition}

The situation is described in the following diagram, where the 
family $g_{\alpha}$ is compatible with the family $f_{\alpha}$:
$$        
\xymatrix@1
        {
         X_{\alpha}\;\; \ar @<+2pt> `u[r] `[rr]^{g_{\alpha}} [rr]
                                             \ar[r]^{f_{\alpha}} 
         & \;\;X\;\;  \ar@{-->}[r]^{\exists ! g} 
         & \;\;Y      
        }
$$

\vspace{1ex}

\begin{definition} \label{terminology} \label{uniqueuptoiso}
  \label{creation} 
Given a functor
$u:\mathbb{T}\to\mathbb{S}$:
\begin{enumerate}
\item  
An object $X$ in $\mathbb{T}$ \emph{sits over} an object $S$ in 
$\mathbb{S}$ when $u(X)=S$. We say also that $X$ is an object 
\emph{over} $S$. An arrow $f:X\to Y$ in $\mathbb{T}$ \emph{sits over}
an arrow $\varphi:S\to T$ in $\mathbb{S}$ when $u(f)=\varphi$, so  
that $X$ (resp. $Y$) sits over $S$ (resp. $T$). We say also that 
$\varphi$ \emph{lifts to} an arrow in  $\mathbb{T}$ when there exists
$f$ over $\varphi$. A family $f_\alpha:X\to Y$ in $\mathbb{T}$  
\emph{sits over} a family $\varphi_\alpha:S\to T$ in $\mathbb{S}$ 
when $u(f_\alpha)=\varphi_\alpha$, for any $\alpha$. We say also that
the family $\varphi_\alpha$ \emph{lifts to} a family in $\mathbb{T}$
when there exists $f_\alpha$ over  $\varphi_\alpha$. 

\item We say that two families $f_{\alpha}:X_{\alpha}\to X$, 
$g_{\alpha}:X_{\alpha}\to Y$ in $\mathbb{T}$ which sit 
over the same family $\phi_{\alpha}:S_{\alpha}\to S$ in $\mathbb{S}$  
are u-isomorphic if there exists an isomorphism $\theta:X\to 
Y$ over $id:S\to S$ such that $\theta\circ f_{\alpha}=g_{\alpha}$ for 
all $\alpha$.
 
\item 
Let $\cc{A}$ be a collection of classes of families in 
$\mathbb{T}$. We say that \mbox{$\cc{A}$-families} are unique up to
isomorphisms if given any two $\cc{A}$-families 
\mbox{$f_{\alpha}:X_{\alpha}\to X$,} $g_{\alpha}:X_{\alpha}\to Y$
which sit over the same family $\phi_{\alpha}:S_{\alpha}\to S$ in
$\mathbb{S}$, they are \mbox{u-isomorphic} by a unique isomorphism. 
(Notice that when $u$ is faithful strict epimorphic families are
unique up to \mbox{isomorphisms.}) 

\item
Consider collections $\cc{A}$ in $\mathbb{T}$ and 
$\cc{B}$ in $\mathbb{S}$. We say that $u$ creates 
$\cc{A}$-families over $\cc{B}$-families if given any 
$\cc{B}$-family $\phi_{\alpha}:S_{\alpha}\to S$, and an object 
$X_{\alpha}$ over $S_{\alpha}$ for every $\alpha$, there exists an 
$\cc{A}$-family $f_{\alpha}:X_{\alpha}\to X$ over 
$\phi_{\alpha}:S_{\alpha}\to S$.

When the class $\cc{B}$ in $\bb{S}$ is the ``same'' class than the
class $\cc{A}$ in $\bb{T}$ (that is, if they are denoted by the same
letter), we simply say that $u$ creates $\cc{A}$-families. 
\end{enumerate}
\end{definition}

All collections considered in this paper are assumed to be closed
under \mbox{u-isomorphisms} without need to say so explicitly.

\begin{definition}\label{surjective} \label{final}
We consider notions for families in  $\mathbb{T}$ relative to the
functor $u$: 
\begin{enumerate} 
\item 
Given a functor $u:\mathbb{T}\to\mathbb{S}$, we
say that a family $f_{\alpha}:X_{\alpha}\to X$ in $\mathbb{T}$ is  
u-surjective when the family $u(f_{\alpha})$ is strict 
epimorphic in $\mathbb{S}$.

\item
Given a functor $u:\mathbb{T}\to\mathbb{S}$, let 
$f_{\alpha}:X_{\alpha}\to X$ be a family in  $\mathbb{T}$ over 
$\varphi_{\alpha}:S_{\alpha}\to S$ in $\mathbb{S}$. The family 
$f_{\alpha}$ is u-final if  for any family 
$g_{\alpha}:X_{\alpha}\to Y$ in $\mathbb{T}$ and arrow $\phi:S\to 
T$ in $\mathbb{S}$  such that $g_\alpha$ sits over 
$\phi\circ\varphi_\alpha$, there exits a unique $g:X\to Y$ over 
$\phi$ such that $g\circ f_\alpha=g_\alpha$.
(final families are unique up to isomorphisms in the sense of
definition \ref{uniqueuptoiso}) 

\item By $\cc{F}\cc{S}=\cc{F} \cap \cc{S}$ we shall denote the
collection of all final and surjective families. 

\end{enumerate}
\end{definition}
We shall often omit the $u$ when we write ``$u$-surjective'' or ``$u
$-final''.
The situation for final families is described in the
following double diagram, where the top diagram sits over the bottom
diagram. 
\begin{equation}   \label{doublediagram}     
\xymatrix@1@R=15pt
        {
         X_{\alpha}\;\; \ar @<+2pt> `u[r] `[rr]^{g_{\alpha}} [rr]
                                             \ar[r]^{f_{\alpha}} 
         & \;\;X\;\;  \ar@{-->}[r]^{\exists ! g} 
         & \;\;Y  
         \\        
         S_{\alpha}\;\; \ar[r]^{\varphi_{\alpha}} 
         & \;\;S\;\;  \ar[r]^{\phi} 
         & \;\;T      
        }
\end{equation}

\vspace{1ex}

Let $s\cc{E}_X$ be the class of all strict 
epimorphic families with codomain $X$. The collection 
$s\cc{E}$ satisfies (I) and (S), but in general it fails to satisfy 
(C), (U) and (F).

Let $\cc{S}_X$ be the class of all surjective 
families with codomain $X$. The collection \mbox{$\cc{S}$} 
satisfies conditions (I) and (S), but fails in general
to satisfy \mbox{(C), (U) and (F).} 

Let $\cc{F}_X$ be the class of all final 
families with codomain $X$. The collection $\cc{F}$ 
satisfies conditions (I), (C) and (S), but fail in 
general to satisfy (U) and (F). 

Finally, the collection $\cc{FS}$ 
satisfies conditions (I) and (S) by \ref{collectionfacts} (\ref{intersection}). 

\vspace{1ex}

\vspace{1ex}

We are interested in the equivalence  
$$\;\;Strict \; epimorphic \;\; \iff \;\; Final\;
surjective\,.$$ Concerning this we have:


 
\begin{fact} \label{characterization1}
If the functor $u$ is faithful, and has left and right adjoints
(notice that this hypothesis is self-dual), then:

A family is strict epimorphic if and only if is final  
and surjective, that is, \mbox{$s\cc{E}=\cc{FS}$.} And the dual statement:
A family is strict monomorphic if and only if is initial  
and injective.
\end{fact}

\begin{definition} \label{efunctor} \label{fmfunctor}
A functor $\bb{T} \mr{u} \bb{S}$ is a $\cc{E}$-functor (resp.
$\cc{M}$-functor) if it is \mbox{faithful} and creates and preserves strict
epimorphic families (resp. strict monomorphic families). If we consider only finite families, we
have the notions of  \mbox{$\cc{E}_{fin}$-functor} and
$\cc{M}_{fin}$-functor. 
\end{definition}

\begin{fact}[about $\cc{E}$-functors] \label{characterization2}
Given a $\cc{E}$-functor (resp. $\cc{E}_{fin}$-functor)  $\bb{T} \mr{u} \bb{S}$, a
family (resp. finite family) in $\bb{T}$
is strict epimorphic if and only 
if is final  and surjective. 
\end{fact}

\begin{fact} [about $\cc{M}$-functors]\label{finitecharacterization2} 
Given a $\cc{M}$-functor (resp. $\cc{M}_{fin}$-functor)  $\bb{T} \mr{u}
\bb{S}$, a family (resp. finite family) in
$\bb{T}$ is strict monomorphic if and only 
if is initial and injective.
\end{fact}   

\vspace{1ex}

In \cite{DE} we have extensively developed a notion of topological
functor. We recall here a couple of results we shall
need in this paper.

\begin{facts}[about topological functors]  \label{topological}
The following holds:
\begin{enumerate}
\item \label{topin}
A functor $u:\mathbb{T}\to\mathbb{S}$ is  
topological if and only if it creates initial families. 
\item \label{topselfdual}
A functor $u:\mathbb{T}\to\mathbb{S}$ is topological if and only if
considered as a functor $u:\mathbb{T}^{op} \to\mathbb{S}^{op}$ is
topological.
\item \label{chartop3}
A functor $u:\mathbb{T}\to\mathbb{S}$ is topological if and only if:

i) It is faithful and preserves and creates strict epimorphic
  families.

ii) It has a full and faithful left adjoint $(-)_\bot \dashv
\, u$, $u(-)_\bot = id$.

\end{enumerate}
\qed \end{facts} 
Topological functors are, in particular, $\cc{E}$-functors and  $\cc{M}$-functors.

\section{$\cc{E}$-functors}  \label{Efunctors}
\vspace{1ex}

The notions of $\cc{E}$-functor and  $\cc{M}_{fin}$-{functor}, $\bb{T}
\mr{u} \bb{S}$, determine a framework of 
the right generality for several constructions found in many
particular situations. Recall that in this case, by
\ref{characterization2} and \ref{finitecharacterization2}, in the
category $\mathbb{T}$ 
strict epimorphic families are the same that final surjective
families, and finite strict monomorphic families are the same that
finite initial injective families. 

An $\cc{E}$-functor does not necesarially \emph{reflect} the property
of being a strict epimorphic family (example, the forgetful functor
of the category of all topological spaces). We establish now some
technical results for future use.

$\cc{E}$-functors (resp $\cc{M}_{fin}$-functors) create any
colimit (resp. finite limit) that may exists in $\bb{S}$, and preserve
any colimit (resp. finite limit) that may exists in $\bb{T}$,
\emph{provided} the colimit (resp. finite limit) of the underlying
diagram already \mbox{exists in $\bb{S}$.}

\begin{proposition} \label{preservescolandlim}
Let $\bb{T} \mr{u} \bb{S}$ be a $\cc{E}$-functor
  (resp. $\cc{M}_{fin}$-functor). Consider a diagram (resp. finite
  diagram) $\Gamma \mr{X} \bb{T}, \; \alpha \mapsto X_\alpha$ in 
$\bb{T}$, and the underlying diagram $S_\alpha = uX_\alpha$ in
  $\bb{S}$. Suppose a colimit cone $S_\alpha \to L$ (resp. a limit
  cone $L \to S_\alpha$) exists in $\bb{S}$. Then: 

i) A colimit cone $X_\alpha \to Y$ over $S_\alpha \to L$ (resp. a
limit cone $Y \to X_\alpha$ over  $L \to S_\alpha$) exists in $\bb{T}$.

ii) If a colimit cone  $X_\alpha \to Y$  (resp. a limit cone $Y \to
X_\alpha$) exists in $\bb{T}$, then $u(X_\alpha \to Y) = S_\alpha \to
L$, (resp. $u(Y \to X_\alpha) = L \to S_\alpha$).
\end{proposition}
\begin{proof}
The reader should be able to carefully check the validity of the two
statements. The faithfulness of $u$ as well as facts
\ref{characterization2} (resp. \ref{finitecharacterization2}) are needed. 
\end{proof}

The proof of the  following proposition is immediate.
\begin{proposition} \label{efunctoryoga1}
Given functors 
$\mathbb{T} \mr{u} \mathbb{S} \mr{v} \mathbb{R}$, if $u$ and $v$ preserve
\mbox{$s\cc{E}$-families} (resp.  $s\cc{M}_{fin}$-families), and the
composite $v \circ u$ creates $s\cc{E}$-families (resp.
\mbox{$s\cc{M}_{fin}$-families),} then $u$ is a $\cc{E}$-functor
(resp. $\cc{M}_{fin}$-functor). 
\qed \end{proposition}

\begin{proposition} \label{efunctoryoga2}
Consider $u:\mathbb{T}\to\mathbb{S}$, $X$ an object in $\mathbb{T}$ 
and the induced functor
$u_{_{\mbox{$*$}}}:\mathbb{T}/_{\mbox{$X$}}
\to   \mathbb{S}/_{\mbox{$S$}}$, where  
$S=u(X)$. Then, if $u$ is an $\cc{E}$-functor so it is
$u_{_{\mbox{$*$}}}$. 
\end{proposition}  
\begin{proof} 
It is clear that $u_{_{\mbox{$*$}}}$ is faithful and that it 
preserves $s\cc{E}$-families. Now we prove that $u_{_{\mbox{$*$}}}$
creates $s\cc{E}$-families. Let $f_\alpha:X_\alpha\to X$ be a family
in $\mathbb{T}/_{\mbox{X}}$ over $\sigma_\alpha:S_\alpha\to S$ in
$\mathbb{S}/_{\mbox{$S$}}$, and an $s\cc{E}$-family $\varphi_\alpha$
in $\mathbb{S}/_{\mbox{$S$}}$ as in the left diagram below:
$$
\xymatrix@R=15pt@C=7pt
        {
         S_\alpha \ar[rr]^{\varphi_\alpha} \ar[dr]_{\sigma_\alpha} 
        &&  R \ar[dl]^{\rho}
        \\
        & S 
        }
\hspace{12ex}
\xymatrix@R=15pt@C=7pt
        {
         X_\alpha \ar@{-->}[rr]^{g_\alpha} \ar[dr]_{f_\alpha} 
        &&  Y \ar@{-->}[dl]^{g}
        \\
        & X 
        }
$$
The canonical functor $\mathbb{S}/_{\mbox{$S$}} \to \mathbb{S}$
preserves $s\cc{E}$-families, so that $\varphi_\alpha$ is 
an $s\cc{E}$-family in $\mathbb{S}$. Since $u$ creates
$s\cc{E}$-families, we have an $s\cc{E}$-family $g_\alpha:X_\alpha\to
Y$ in $\mathbb{T}$ over $\varphi_\alpha:S_\alpha\to R$. 
Moreover, since the family $g_\alpha$ is final surjective in
$\mathbb{T}$, there exists a unique $g:Y\to X$ such that $g\circ  
g_\alpha=f_\alpha$ for all $\alpha$. It  
is immediate to check that the family $g_\alpha$ in 
$\mathbb{T}/_{\mbox{$X$}}$ over $\varphi_\alpha$ in
$\mathbb{S}/_{\mbox{$S$}}$ is an  $s\cc{E}$-family in 
$\mathbb{T}/_{\mbox{$X$}}$.     
\end{proof}

\begin{definition} \label{ffactorizations} 
A category $\mathbb{T}$ has $f$-factorizations if every family
 $g_\alpha:X_\alpha\to X$ factorizes in the form $g_\alpha=m\circ
 h_\alpha$,
 where $m:H\to X$ is a monomorphism and $h_\alpha:X_\alpha\to H$ is an 
strict epimorphic family.
\end{definition}

\begin{remark} \label{miso} 
If a family $g_\alpha:X_\alpha\to X$ factorizes in the form 
$g_\alpha=m\circ h_\alpha$, where $m:H\to X$ is a monomorphism, then a 
family  $f_\alpha:X_\alpha\to Y$ is compatible with $g_\alpha$ if and
only 
if it is compatible with $h_\alpha$. Moreover, if $g_\alpha=m\circ
h_\alpha$ 
is a $f$-factorization then $g_\alpha$ is a strict  epimorphic family
if 
and only if $m$ is an isomorphism.   
\qed\end{remark}

\begin{proposition} \label{factorizationimpliesC}
If a category $\mathbb{T}$ has $f$-factorizations, then strict
epimorphic families compose, that is, the collection $s\cc{E}$ has
property (C).
\end{proposition}
 
\begin{proof} 
We consider a family $X_\alpha \mr{ g_\alpha} X$ and, for any
$\alpha$, 
a family $X_{\alpha,\beta} \mr{g_{\alpha,\beta}} X_\alpha$, so that we
have 
the composite family $l_{\alpha,\beta}=g_\alpha\circ g_{\alpha,\beta}:
X_{\alpha,\beta} \to X$. Let us suppose that $g_\alpha$ and
$g_{\alpha,\beta}$ are strict epimorphic families. To prove that so is
$l_{\alpha,\beta}$ we 
take the $f$-factorization $l_{\alpha,\beta}=m\circ h_{\alpha,\beta}$, 
with $H\mr{m} X$ mono and $h_{\alpha,\beta}$ an strict epimorphic
family. 
Fixing  $\alpha$, it is clear that the family $h_{\alpha,\beta}$ is
compatible with $g_{\alpha,\beta}$, hence there exists an arrow
$X_\alpha \mr{h_\alpha} H$ unique such that $h_\alpha\circ
g_{\alpha,\beta}=h_{\alpha,\beta}$. Moreover (compose with the
epimorphic family 
$g_{\alpha,\beta}$) we have $m\circ h_\alpha=g_\alpha$ for any
$\alpha$. Now 
we prove that the family $h_\alpha$ is strict epimorphic family. In
fact: if a family $f_\alpha$ is compatible with $h_\alpha$ then the 
family $f_{\alpha,\beta}=f_\alpha\circ g_{\alpha,\beta}$ is compatible 
with $l_{\alpha,\beta}$, so that there exists a unique arrow $f_H$
such 
that $f_H\circ h_{\alpha,\beta}=f_{\alpha,\beta}$, hence (compose with
the epimorphic family $g_{\alpha,\beta}$) $f_H\circ
h_\alpha=g_\alpha$. 
$f_H$ is unique with this condition because $h_\alpha$ is an
epimorphic family (notice that so is $h_{\alpha,\beta}$). Finally,
remark \ref{miso} implies that $m$ is an isomorphism, hence
$l_{\alpha,\beta}$ is an strict epimorphic family. 
\end{proof}

\begin{proposition} \label{liftfactorizations}
If $u$ is an $\cc{E}$-functor and $\mathbb{S}$ has 
$f$-factorizations (definition \ref{ffactorizations}), then $\bb{T}$
has $f$-factorizations.  
\end{proposition}
\begin{proof}
 Given any family $g_\alpha:X_\alpha \to X$ in $\mathbb{T}$, we can
 take the $f$-factorization \mbox{$u(g_\alpha)=m\circ h_\alpha:u(X_\alpha)
 \to H\to u(X)$} in $\mathbb{S}$, and a strict epimorphic family
 \mbox{$f_\alpha:X_\alpha \to Y$} in $\mathbb{T}$ over $h_\alpha$. But
 $h_\alpha$ is final surjective, so that there exists $f$ over $m$
 such that $g_\alpha=f\circ f_\alpha$, with $f$ mono because $u$ is
 faithful.
\end{proof}                        

It is well known and easy to prove that if a functor has a right
adjoint then it preserves strict epimorphic families. Next we
establish that under certain hypotheses this condition is also
sufficient for the existence of a right adjoint. Moreover, we give an
explicit construction of this adjoint. It is convenient to 
consider before some particular cases.

\begin{proposition} \label{preservesradjoint1}
If $u$ is an $\cc{E}$-functor and $\mathbb{S}$ has 
$f$-factorizations (definition \ref{ffactorizations}), then $u$ has a
right adjoint. 
\qed\end{proposition}
We indicate the idea of the proof. The construction of the right
adjoint is as follows: Given an  
object $S$ in $\mathbb{S}$, consider in $\mathbb{S}$ all 
arrows of the form $u(X) \mr{\varphi} S$. 
Factor this family 
\mbox{$u(X) \mr{\psi} H \mr{m}  S$,} with a monomorphism 
$m$ and the family $\psi$ strict epimorphic. Let
$R(S)$ 
and $X \mr{\bar{\varphi}} R(S)$ be a strict epimorphic
  family over the family
$\psi$.  Then, R(S) is the right adjoint to $u$ on $S$.

\begin{theorem}\label{preservesradjoint2}
Consider a diagram of categories and functors where $u$ and $u'$ are 
$\cc{E}$-functors, $u'\circ F=u$.
$$
\xymatrix@R=15pt@C=7pt
        {
         \mathbb{T} \ar[rr]^{F} \ar[dr]_{u} 
        && \mathbb{T}' \ar[dl]^{u'}
        \\
        & \mathbb{S} 
        }
$$
Assume that 
$\mathbb{S}$ has $f$-factorizations (definition
\ref{ffactorizations}). Then $F$ has a right  
adjoint if and only if $F$ preserves $s\cc{E}$-families.  
\qed\end{theorem}
We indicate the idea of the proof. The construction of the right
adjoint is as follows: Given an object 
$Y$ in $\mathbb{T}'$ consider 
in $\mathbb{S}$ all arrows of the form \mbox{$u(X) \mr{\varphi} u'(Y)$} such
that $\varphi=u'(g)$ for some $F(X) \mr{g} Y$. Factor this family 
\mbox{$u(X) \mr{\psi} H \mr{m}  u'(Y)$,}
 for a monomorphism $m$ and an
strict epimorphic family  $\psi$. Let 
$G(Y)$ and  $X \mr{\bar{\varphi}} G(Y)$ be an strict epimorphic 
family 
over  the family $\psi$. Then, $G(Y)$ is the right adjoint to 
\mbox{$F$ on $Y$.} 

\vspace{1ex}

Now we set and prove in detail the general theorem: 

\begin{theorem}\label{pra3}
Consider a diagram of categories and functors where $u$ and $u'$ are 
$\cc{E}$-functors, $u'\circ F=L\circ u$, and $R\vdash L$. 
$$
\xymatrix
        {
         \mathbb{T} \ar[r]^{F} \ar[d]_{u} 
        & \mathbb{T}' \ar[d]^{u'}
        \\
        \mathbb{S} \ar@<-1ex>[r]_{L} 
        & \mathbb{S}' \ar@<-1ex>[l]_{R} 
        }
$$
Assume that
$\mathbb{S}$ has $f$-factorizations (definition 
\ref{ffactorizations}). Then $F$ has a right 
adjoint if and only if $F$ preserves $s\cc{E}$-families. 
\end{theorem}
\begin{proof} 
Suppose that $F$ preserves strict epimorphic families. 
For every object $Y$ in $\mathbb{T}'$ we construct an object $G(Y)$ 
in  $\mathbb{T}$ and an arrow $\epsilon_Y:F(G(Y))\to Y$ universal 
from $F$ to $Y$. First, consider  arrows  
$\varphi:u(X)\to R(u'(Y))$ in $\mathbb{S}$. To every $\varphi$ it
corresponds 
under the adjunction $R\vdash L$ an arrow
\mbox{$\hat{\varphi}:L(u(X))\to u'(Y)$.}  Now, we take in 
$\mathbb{S}$ the family of all arrows $\varphi$ as above such
that $\hat{\varphi}=u'(g)$ for some $g:F(X)\to Y$ in
$\bb{T}'$. Notice that we need the equation $u'\circ F=L\circ
u$. Since 
$\mathbb{S}$ has $f$-factorizations, there exists a
monomorphism $m:H\mono R(u'(Y))$ such that 
$\varphi=m\circ\psi\,:\, u(X) \mr{\psi} H \mr{m} R(u'(Y))$ for
all $\varphi$ in our family, and the family $\psi$ is strict 
epimorphic. But $u$ creates $s\cc{E}$-families, hence there exists an
object $G(Y)$ in $\mathbb{T}$ and for any $\varphi$ in the family an 
arrow  $\bar{\varphi}:X\to G(Y)$ such that $\bar{\varphi}$ is an 
strict epimorphic family over the family $\psi$ (in particular 
$u(G(Y))=H$). By hypothesis, $F(\bar{\varphi}):F(X)\to F(G(Y))$ is 
a \mbox{$s\cc{E}$-family} in $\mathbb{T}'$, and it sits over the
family $L(\psi)$. Moreover it is final and surjective because
$u'$ is an $\cc{E}$-functor. Recall that the family $g:F(X)\to Y$ 
considered above sits over the family $\hat{\varphi}$ which 
factorizes trough $L(\psi)$, so that there exists a unique 
$\epsilon_Y:F(G(Y))\to Y$ such that it sits over $\hat{m}$ and 
$\epsilon_Y\circ F(\bar{\varphi})=g$. It remains to prove that this 
arrow is universal. The above construction shows that given $g$ we 
have $\bar{\varphi}$ such that $\epsilon_Y\circ
F(\bar{\varphi})=g$. Consider any other arrow $f$ such that 
$\epsilon_Y\circ F(f)=g$.  Then $\hat{m}\circ L(u(f))=u'(g)$, and 
applying the adjunction it follows $m\circ u(f)=\varphi=m\circ 
u(\bar{\varphi})$. Thus, $f=\bar{\varphi}$ because $m$ is a
monomorphism and $u$ is faithful.         
\end{proof}
A functor $F$ as above which has a right adjoint in particular
preserves strict epimorphic families, thus it will preserve final
surjective families (which are the same). However, \emph{it will not
  preserve arbitrary final families in general}.

A first application of theorem \ref{pra3} is the
following: 
\begin{corollary} \label{corollarypra3} 
Let $\mathbb{S},\mathbb{T}$ be categories with finite products and 
$u:\mathbb{T}\to\mathbb{S}$ an \mbox{$\cc{E}$-functor} which 
preserves   
finite products. Assume that $\mathbb{S}$ is cartesian closed and
that it  
has $f$-factorizations. Then, $\mathbb{T}$ is cartesian closed 
if and only if the cartesian product in $\mathbb{T}$ preserves 
$s\cc{E}$-families.
\end{corollary}
\begin{proof} 
Apply Theorem \ref{pra3} with the following diagram of functors:
$$
\xymatrix
        {
         \mathbb{T} \ar[r]^{(-)\times X} \ar[d]_{u} 
        & \mathbb{T} \ar[d]^{u}
        \\  
        \mathbb{S} \ar@<-1ex>[r]_{(-)\times S}  
        & \mathbb{S} \ar@<-1ex>[l]_{(-)^S} 
                     \ar@{}[r] 
        & (where \; u(X)=S).
        }
$$ 
\end{proof} 
It is interesting to give some further details of the construction of
the exponential in $\bb{T}$  that follow from the proof of theorem
\ref{pra3}. Given any two objects $X$ and $Y$ in $\bb{T}$, 
the underline object in $\bb{S}$ of the exponential $Y^X$ is a
subobject of the exponential between the underline objects,  $u(Y^X)
\mono u(Y)^{u(X)}$. The family of all arrows
$Z \mr{\bar{\varphi}} Y^X$  
such that the arrow $u(Z)\times u(X)\to u(Y)$ (which corresponds by
adjointness  to
the arrow $u(Z) \mr{u(\bar{\varphi})} u(Y^X) \mono u(Y)^{u(X)}$) 
lifts into $Z\times X\to Y$, is a final surjective family.

We now  generalize this corollary to
localized exponentials. Recall  that if $\mathbb{S}$ has finite
products,  then the usual functor $\mathbb{S}/S \to \mathbb{S}$ has a right adjoint, 
and that $\mathbb{S}$ has finite limits if and only if each 
$\mathbb{S}/S$ does. Recall also that by definition a category
$\mathbb{S}$ is  
\emph{locally cartesian closed} if $\mathbb{S}/S$ is cartesian 
closed for any object $S$ in $\mathbb{S}$. Given an arrow $R \to S$
 in a locally cartesian  
closed category $\bb{S}$, the pulling-back functor 
$\mathbb{S}/S\to\mathbb{S}/R$ has a right adjoint.  

\begin{theorem} \label{lcc}
Let $\mathbb{S},\mathbb{T}$ be categories with finite limits and 
$u:\mathbb{T}\to\mathbb{S}$ an \mbox{$\cc{E}$-functor} which preserves 
finite limits. Assume that $\mathbb{S}$ is locally cartesian closed 
and has \mbox{$f$-factorizations.} Then, $\mathbb{T}$ is locally 
cartesian closed if and only if $s\cc{E}$-families are universal in 
$\mathbb{T}$.
\end{theorem}
\begin{proof} 
Notice that if $\mathbb{S}$ has $f$-factorizations then so does 
each localized category $\mathbb{S}/S$. Then apply corollary
\ref{corollarypra3} and proposition \ref{efunctoryoga2}.
\end{proof} 

Our last result in this section concerns subobject classifiers.

\begin{proposition} \label{liftomega}
Let $\mathbb{S},\mathbb{T}$ be categories with finite limits and 
$u:\mathbb{T}\to\mathbb{S}$ a \mbox{$\cc{M}_{fin}$-functor} with a right adjoint
$R$. Assume that $\mathbb{S}$ has $f$-factorizations and a strict subobject 
classifier $\Omega$. Then, $R\Omega$ is an strict subobject classifier
in $\bb{T}$.
 \end{proposition}
\begin{proof}
Let $1 \mr{t} \Omega$ be the generic strict subobject. Clearly $1 =
  R1$, and there is a map $1 \mr{at} R\Omega$ in $\bb{T}$
  corresponding by adjointness to the map 
$1 \mr{t} \Omega$ in $\bb{S}$ (notice that since $u$ preserves finite
limits (proposition \ref{preservescolandlim}), $u1 = 1$). Given a strict subobject
  $M \mono X$ in $\bb{T}$, the strict subobject  $uM \mono uX$ in
  $\bb{S}$ determines the pullback square on the left below, which in
  turn determines by adjointness the commutative square on the right.
$$
\xymatrix@R=15pt@C=20pt
              {
                   uM \ar[r]  \ar[d]  
               &  uX  \ar[d]^\varphi
               \\
                  1 \ar[r]^{t}  
               & \Omega
              }
\hspace{12ex}
\xymatrix@R=15pt@C=20pt
              {
                   M \ar[r]  \ar[d]  
               &  X  \ar[d]^{a\varphi}
               \\
                  1 \ar[r]^{at}  
               & R\Omega
              }
$$   
It remains to see that this square is a pullback. This follows by the
fact that the strict subobject $M \mono X$ is initial injective
(\ref{finitecharacterization2}).
\end{proof}

\section{$f$-regular categories and $f$-quasitopoi}  \vspace{1ex}   \label{f-qandqfunc}

Recall that a  \emph{regular category} is a category with finite
limits and such that any arrow can be factorized into a monomorphism
composed with an  strict epimorphism, and in addition strict
epimorphisms  are universal. In a regular category strict epimorphisms
are  the same that regular epimorphisms, and strict epimorphisms
compose  \cite{B}, \cite{G}. We introduce now a notion which corresponds to
the  notion of regular category, but utilizing strict epimorphic
families  instead of single strict epimorphisms. This notion is not
elementary and it means a (co) completeness requirement. We call this notion
\emph{f-regular}, \emph{f} for family.  

\begin{definition} \label{f-reg}
A category is f-\emph{regular} when  it satisfies:

R1) It has all finite limits.

R2) Strict epimorphic families  are 
universal (U in definition \ref{familyproperties}).

R3) It has $f$-factorizations (definition \ref{ffactorizations}).

\end{definition} 

Notice that (see definition \ref{familyproperties}) the second
condition in the definition of \mbox{$f$-regular} category means that
the class $s\cc{E}$ of strict epimorphic families is stable under
pulling-back. In the following remark we set an essential property of
\mbox{$f$-regular categories.} It follows 
from \ref{collectionfacts} and
\ref{factorizationimpliesC}: 
 
\begin{remark}\label{goodfamily}
In a $f$-regular category the collection $s\cc{E}$ of strict epimorphic families satisfies all five properties 
in definition \ref{familyproperties}.  
\qed \end{remark}

From proposition \ref{liftfactorizations} we have:

\begin{proposition} \label{liftfregular}
Given a $\cc{E}$-functor $\mathbb{Q}\mr{u}\mathbb{S}$, with
$\bb{S}$ a 
$f$-regular category, then $\bb{Q}$ is also a $f$-regular category if
and only if the families created by $u$ are universal
\qed\end{proposition}

Clearly, $f$-regular categories are regular.  Regular categories
sufficiently cocomplete (for example if the lattice of subobjects have
arbitrary suprema) are \mbox{$f$-regular.} However, we consider
$f$-regularity to be the primitive  notion since when working with
families this is the notion which arises naturally.

\begin{proposition} \label{fregular2}
A category is $f$-regular if and only if it satisfies:

R1) It has all finite limits.

R2) Strict epimorphic families  are 
universal (U in definition \ref{familyproperties}).

R4) It has coequalizers of kernel pairs.

R5) The lattice of subobjects of any object has arbitrary suprema.
\end{proposition}
\begin{proof}
Under the presence of R1) and R2) we have:

R3) $\Rightarrow$ R4): Given an arrow $S \mr{f} T$, we can factorize
it as \mbox{$S \mr{g} H \mono T$,} with $g$ a strict
epimorphism. Then, it readily follows that $S \mr{g} H$ is a
coequalizer of the kernel pair of $f$. 

R3) $\Rightarrow$ R5):  Given a family of subobjects $m_\alpha:H_\alpha \mono S$, we can factorize it as
$H_\alpha \mr{i_\alpha} H \mono S$, with
$i_\alpha$ an strict epimorphic family (note that each $i_\alpha$ is
mono). Then, it readily follows that $H=\bigvee_\alpha H_\alpha$.

R4), R5)  $\Rightarrow$ R3). Given a family $X_\alpha \mr{g_\alpha} X$,
take for any $\alpha$ the coequalizer of the kernel pair of $g_\alpha$.
We have a factorization 
$X_\alpha \mr{h_\alpha} H_\alpha \mr{m_\alpha} X$, with $h_\alpha$ a
strict epimorphism. It is known that
from R2) it follows that  $m_\alpha$ is a mono (this is
an argument due originally to M. Tierney and independently M. Kelly
\cite{MK}, see  2.1.2 and 2.1.3 in \cite{BO}). Let $H =
\bigvee_\alpha H_\alpha$, abuse notation, and write $X_\alpha
\mr{h_\alpha} H$. Then, it can be proved that $h_\alpha$ is a strict
epimorphic family. 

\end{proof} 

\vspace{1ex}

Given a pointed (small) site, the category of separated sheaves is a
\emph{bounded quasitopos} in the sense of Penon \cite{P2}. In
practice the most conspicuous quasitopoi are not bounded. They are
categories of separated 
sheaves for pointed \emph{large} sites. In this case, the categories
of sheaves are not topoi, and they are not even legitimate categories
since the class of morphisms between two sheaves is not in general a
set. However, separate sheaves form legitimate categories 
(insofar the
hom-sets are small) which are elementary quasitopoi in the sense of
Penon. They bear to elementary quasitopoi a relation which should
be considered as corresponding (in the context of quasitopoi) to the
relation that Grothendieck topoi bear to elementary 
topoi. In this context, the size condition (``bounded'') is
unnecessary, and probably misleading. Furthermore, it is not satisfied
by many important examples. We introduce now the
abstract notion which describes this situation.

\begin{definition} \label{f-q}
A category $\mathbb{Q}$ is an \emph{f-quasitopos} if it satisfies:

(QT1) It has all finite limits.

(QT2) It has all small colimits.

(QT3) It is locally cartesian closed.

(QT4) It has an strict subobject classifier $1 \mr{t} \Omega$.

(QT5) It has f-factorizations.  
\end{definition}

Since (QT3) implies that strict epimorphic families are universal, 
$f$-quasitopoi are, in particular, $f$-regular categories.

Recall that a Penon's (or \emph{elementary}) quasitopos \cite{P2}, is a category
which satisfies (QT1), (QT3), (QT4),  and a elementary (and weaker)
form of (QT2), namely: it has all finite colimits. A
$f$-quasitopos is, in particular, an elementary quasitopos.

For a \emph{bounded} quasitopos Penon requires QT1 to QT4, and a
different (and stronger) form of (QT5), namely: for any object $S$,
the lattice of 
subobjects $P(S)$ should be small. This condition implies that $P(S)$
is a complete lattice by the existence of small colimits. From
proposition \ref{fregular2} it follows:   

\begin{proposition} \label{onQT5}
In definition \ref{f-q}, condition (QT5) is equivalent to: For any object $S$,
the lattice $P(S)$ of all subobjects of $S$ is a complete lattice.   
\qed\end{proposition}
 
Remark that the existence of an initial object by (QT2) furnishes the
factorization of the empty family in (QT5), and vice-versa.  

\vspace{1ex}

Between the three notions of quasitopoi, we have the
(strict) implications:  

$$\text{bounded quasitopos} \;\;\rimply \;\;f\text{-quasitopos} \;\;\rimply\;\;
\text{elementary quasitopos}.$$   
 
\vspace{1ex}

Now we consider a kind of functor $u:\mathbb{Q}\to\mathbb{S}$ that
correspond to inverse images of geometric morphisms of topoi. These
functors are relevant when the categories are $f$-regular. 

\begin{definition} \label{fregularfunctor1}
A $\cc{M}_{fin}$ and $\cc{E}$-functor $u:\mathbb{Q}\to\mathbb{S}$ between
$f$-regular categories is called a \emph{$f$-regular functor}. 
\end{definition}

The $f$-regular functors should be considered as the \emph{morphisms}
of $f$-regular categories and $f$-quasitopoi. The variance in these
large categories should be opposite to the variance in the
category of Grothendieck topoi, since $f$-quasitopoi are categories of
``generalized spaces'', rather than being \mbox{``generalized spaces''
themselves.}

\vspace{1ex} 

Some times it is convenient to consider  $\cc{M}_{fin}$ and
$\cc{E}$-functors between categories which are not necesarially $f$-regular.

\begin{definition} \label{fregularfunctor2}
A $f$-regular functor $u:\mathbb{Q}\to\mathbb{S}$ between
arbitrary categories is a functor such that:

RF1) It has a right adjoint.

RF2) It is faithful.

RF3) It creates and preserves finite strict monomorphic families. 

RF4) It creates and preserves strict epimorphic families. 

RF5) It creates and preserves universal strict epimorphic families. 
\end{definition}

Clearly, when the categories are $f$-regular, both definitions
coincide (recall proposition \ref{preservesradjoint1}).

From proposition \ref{liftfregular} it immediately follows:
\begin{proposition} \label{quasi=fmande}
Given a $f$-regular functor $u:\mathbb{Q}\to\mathbb{S}$, if $\,\bb{S}$
is a $f$-regular \mbox{category,} then so is $\bb{Q}$.
\qed\end{proposition}

\begin{theorem} \label{lift}
Given a $f$-regular functor $u:\mathbb{Q}\to\mathbb{S}$, if
$\,\mathbb{S}$ is is a f-quasitopos then so is $\mathbb{Q}$. 
\end{theorem}
\begin{proof}
Follows by proposition \ref{preservescolandlim}, theorem \ref{lcc} and
proposition \ref{liftomega}.
\end{proof}

The following are basic properties of $f$-regular functors:

\begin{proposition} \label{sm=ii}
Given a $f$-regular functor $\mathbb{Q}\mr{u}\mathbb{S}$ between
arbitrary \mbox{categories,}

1)  A finite family in $\mathbb{Q}$ is strict monomorphic if and only
if it is initial injective.

2)  Any family in $\mathbb{Q}$ is strict epimorphic if and only if it
is final surjective.

3)  Creates any finite limit that may exist in $\bb{S}$, and preserves
    any finite limit that may exist in $\bb{Q}$, provided it already
    exists in $\bb{S}$. 

4) Creates any colimit that may exist in $\bb{S}$, and preserves any
colimit that may exist in $\bb{Q}$. When
the colimit is universal in $\bb{S}$, the created colimit in $\bb{Q}$
is also universal.

5) It has a right adjoint $u \dashv R$, but not necesarially $uR
\cong id$. It does not have in general a left adjoint.

6) It preserves and creates $f$-factorizations. 

7) The usual construction of exponentials in $\bb{Q}$ out of
   exponentials in $\bb{S}$ holds.

8) If $\Omega$ is an strict subobject classifier in $\bb{S}$, then
$R\Omega$ is an strict subobject classifier in $\bb{Q}$.
\qed\end{proposition} 

Topological and $f$-regular functors are a different kind of
$\cc{E}$-functors. A topological functor may fail to satisfy condition
RF5). A $f$-regular functor may not create initial families which are
not finite or injective (that is, over non finite or strict
monomorphic families), or final 
families which are not surjective (that is, over non strict epimorphic
families). However,
often in practice there are functors which 
are topological and $f$-regular simultaneously (for example, the
categories of quasispaces in section \ref{quasispacessection}, or of
strict quasispaces in the context described in \ref{classical}. The following results
clarify this situation. 

\begin{theorem} \label{tisq1}
A topological functor $u:\mathbb{T}\to\mathbb{S}$ satisfies conditions
RF1) to RF4) in definition \ref{fregularfunctor2}. Thus, it is
$f$-regular if (and only if) it satisfies FR5).
\end{theorem}
\begin{proof}
Consider \ref{topological}. By item (\ref{chartop3} i) and its dual
(which holds 
by item (\ref{topselfdual})), a topological functor satisfies RF2), RF3) and
RF4). By  the dual of item (\ref{chartop3} ii), it
satisfies RF1).    
\end{proof}

\begin{corollary} \label{tisq2}
If strict epimorphic families in $\bb{S}$ are universal (in
particular, if $\bb{S}$ is $f$-regular), then: 

A topological functor
$u:\mathbb{T}\to\mathbb{S}$ is $f$-regular if
(and only if) strict epimorphic families in $\mathbb{T}$ are universal. 
\qed\end{corollary}

\begin{theorem} \label{qist}
A $f$-regular functor $u:\mathbb{Q}\to\mathbb{S}$ is topological
if (and only if) it has a left adjoint $(-)_\bot \dashv \, u$, with
$u(-)_\bot = id$. 
\end{theorem}
\begin{proof}
Follows immediately by \ref{topological} $\,$(\ref{chartop3}).
\end{proof}

\section{Quasispaces over a category $\bb{S}$}  \label{quasispacessection}

In this section we develop the constructions in the theory of
quasitopoi which correspond to the construction of categories of
sheaves in the theory of topoi. These constructions, or closely related
ones, were considered by Antoine over the category of Sets \cite{A}, and
by Penon over general categories \cite{P2}. The paradigmatic example
behind all this are Spanier's quasitopologies \cite{S}. 

\vspace{1ex}

\noindent {\bf Notation}
The hom-sets for any category will be denoted with square
brackets. Thus, given any two objects $X, \; Y$ in any category $\bb{X}$, 
$[X,  \, Y] \in \bb{S}et$ denotes the set of arrows in $\bb{X}$ from $X$ to $Y$.
 
\vspace{1ex}

We recall that a Grothendieck pretopology $\cc{J}$ on a category 
$\mathbb{C}$ consists of a family of \emph{covers} $\cc{J}_C$ for each object $C \in \bb{C}$ satisfying properties (I), (C) and (U) in definition \ref{familyproperties}  (thus it also satisfies (F), see \ref{collectionfacts} (\ref{UCimpliesF})). A Grothendieck topology is a pretopology which in addition satisfies property (S). By adding all families which are refined by a cover any pretopology generates a topology with no other additional families.

From now on 
we suppose that the following data are given:
\begin{equation} \label{data}
A \; category \; \mathbb{C} \; 
\; with \; a \; Grothendieck \; topology \; \cc{J},\; and \; a \; functor \; 
u:\mathbb{C}\to\mathbb{S}.
\end{equation}

 For any object $S$ in $\mathbb{S}$, there is 
a presheaf:
$[u(-),\,S]: \mathbb{C}^{op} \longrightarrow  \bb{S}et$, with the usual action on morphisms. For any 
$\varphi:S\to T$ in $\mathbb{S}$, composing with $\varphi$ is a natural 
transformation \mbox{$[u(-),\,S] \mr{\varphi^*} [u(-),\,T]$.}

 \begin{definition} A \emph{quasispace} is a subpresheaf $X \subset [u(-),\,S]$ satisfying a \emph{covering condition}. 
 Given $C \in \bb{C}$, the maps $u(C) \mr{\sigma} S$ in $X(C)$ are called \emph{admissible maps}. We also say that $X$ is a quasispace structure on the set S. 

Covering condition: Given $u(C) \mr{\sigma} S$ and a cover $C_\alpha \mr{f\alpha} C$ in $\cc{J}_C$,  if the composite $\sigma \circ u(f_\alpha )$ is admissible for all $\alpha$, then so is  $\sigma$.
\end{definition} 

It is convenient to lay down this definition explicitly:

\begin{definition} [{\bf explicit}] A \emph{quasispace} is a pair $(S,\, X)$ 
where $S$ is an object of $\mathbb{S}$ and $X$ assigns to each object 
$C \in \bb{C}$ a subset $XC \subset [u(C),\,S]$ subject to the following conditions:

Presheaf condition: 
$$(C \mr{f} D)  \in  \bb{C},\; (uD \mr{\sigma} S) \in  XD  
\;\; \Rightarrow \;\; (uC \mr{uf} uD  \mr{\sigma}  S)  \in XC.
$$

Covering condition: 
$$
(C_\alpha \mr{f_\alpha} C) \in  \cc{J}_C,\; \forall \alpha \: (uC_\alpha \mr{uf_\alpha} uC  \mr{\sigma}  S)  \in XC_\alpha \;\; \Rightarrow \;\;  (uC \mr{\sigma} S) \in  XC.
$$
\end{definition} 

Notice that if the covering condition is satisfied over a pretopology, it will be satisfied also over the generated topology.
 
Notice that the covering condition is not  
the sheaf condition on the presheaf \mbox{$X: \bb{C}^{op} \to \bb{S}$.} Using Grothendieck's abuse of notation $C \mr{\sigma} X$ for \mbox{$(uC \mr{\sigma} S) \in XC$,} the later takes the form:
$$
(C_\alpha \mr{f_\alpha} C) \in  \cc{J}_C,\;  (C_\alpha  \mr{\sigma_\alpha}  X) \; compatible \;\; \Rightarrow \;\;  \exists \,! (C \mr{\sigma} X) \;|\; \sigma \circ f_\alpha = \sigma_\alpha.
$$

 We see then that the sheaf condition will be satisfied by a
 quasispace precisely when the family   
$uC_\alpha \mr{f_\alpha} uC$ is a strict epimorphic family of $\bb{S}$, that is, $C_\alpha \mr{f_\alpha} C$ is a surjective family of $\bb{C}$. We have:
\begin{remark} \label{sheafcondition}
If the covers are surjective families, then the sheaf condition is equivalent to the covering condition. Thus, in this case, a subpresheaf of [u(-),\,S] is a sheaf if and only if it is a quasispace.
\qed\end{remark}

Morphisms of quasispaces are arrows 
 $\varphi$ in $\bb{S}$ such that the natural transformation $\varphi^*$ sends admissible maps to admissible maps. Explicitly:
\begin{definition}
A \emph{morphism of quasispaces} $(S,\, X ) \mr{\varphi} (T,\, Y)$ is an 
arrow \mbox{$S \mr{\varphi} T$}
  in $\mathbb{S}$ such that: 
$$(uC \mr{\sigma} S) \in  XC  
\;\; \Rightarrow \;\; (uC \mr{\sigma} S  \mr{\varphi}  T)  \in YC.
$$ 
\end{definition}

We shall denote $\bb{Q}$ the category of quasispaces. By definition there is a faithful forgetful functor that we denote $\bb{Q} \mr{q} \bb{S}$, $q(S,\,X) = S$, $q(f) = f$.

\vspace{1ex}

On any object $S \in \bb{S}$ there is a maximal and a minimal quasispace structure:

\begin{example} \label{example1}
$
S_\top = (S,\, S_\top), \;\; where  \;  S_\top C = [uC,\,S] .
$

\vspace{1ex} 

$
S_\bot = (S,\, S_\bot), \;\; where \; S_\bot C= 
\left \{
\begin{array}{ll}
[uC,\,S]  & \mbox{if the empty family is in $\cc{J}_C$} \\ 
\emptyset  &  \mbox{otherwise} 
\end{array}
\right.
$
\qed\end{example}

\begin{example} \label{example2}
If $1$ is a terminal object of $\bb{S}$, then $1_\top$ is a terminal object of $\bb{Q}$. Observe that $1_\bot \neq 1_\top$. In general, a quasispace structure $X$ on $1 \in \bb{S}$,
\mbox{$1_\bot  \subset X \subset  1_\top$,} is determined by a \emph{sieve} $\bb{A} \subset \bb{C}$ ($C \in  \bb{A}  \siff  XC = 1$) "closed under covers" (condition which is vacuous when the topology is trivial). We see that in general there is a proper class of different quasispace structures on $1$.
\qed\end{example}

We warn the reader that $1_\bot \neq 1_\top$ even in the case where $\bb{S} = \bb{S}et$ is the category of sets and $1$ is the singleton set. In this case, $1_\bot$ represents the forgetful functor, but it is not the terminal object  $1_\top$ of $\bb{Q}$. Classically a further condition is imposed on a quasispace to have the forgetful functor represented by the terminal object (see definition \ref{strictquasispace}).

\begin{example} \label{zero}
If $0$ is an initial object of $\bb{S}$, then $0_\bot$ is an initial object of $\bb{Q}$. When $0$ is empty (that is, $[S, \, 0]  \neq  \emptyset \siff S = 0)$, we shall denote $0  = \emptyset$. In this case, if in addition, the empty family covers $C$ if and only if $uC = \emptyset$, then,  there is only one quasispace structure on $\emptyset$, and $\emptyset_\bot = \emptyset_\top$ is also an empty initial object of $\bb{Q}$.
\qed\end{example}

\begin{yoneda} $\,$  \label{yoneda} 
\begin{enumerate}
\item \label{y1} 
Given any object $C \in \bb{C}$, $uC$ carries a canonical structure of quasispace, that we denote $\varepsilon C$, defined by stipulating the equivalence:
\begin{center}
\begin{picture}(0, 30) 
\put (-50, 18)  {$uK \mr{\sigma} uC \; \in \; \varepsilon C(K)$}
\put (-145, 12){\line(1, 0){290}}
\put (-130, 0) {$\exists \; K_i \mr{k_i} K \in \, \cc{J}_K \;\;and\;\;  K_i \mr{f_i} C \;\; such \; that\;\; \sigma \circ uk_i = uf_i$}
 \end{picture}
 \end{center}
 
\vspace{1ex}

\noindent (in particular, any $uK \mr{\sigma} uC$ that lifts to $\bb{C}$ is admissible for $\varepsilon C$).   

\vspace{1ex}

\noindent The assignment $C \mapsto (uC, \, \varepsilon C)$ define a
functor (which acts as the equality on arrows) $\bb{C} 
\mr{\varepsilon} \bb{Q}$. We abuse notation and write $\varepsilon C =
(uC, \, \varepsilon C)$. 

\noindent Clearly $q \varepsilon = u$ :  
$
\hspace{12ex}
\xymatrix@R=15pt@C=7pt
         {
          \bb{C} \; \ar[rr]^{\varepsilon}  \ar[rd]_{u}  
           &&  \mathbb{Q} \ar[ld]^{q}  
           \\
           & \bb{S} 
         }
$ 

\vspace{1ex}

\item \label{y2}

\noindent If covers are final families in $\bb{C}$, then

\begin{center}
\begin{picture}(0, 30) 
\put (-50, 18)  {$uK \mr{\sigma} uC \; \in \; \varepsilon C(K)$}
\put (-80, 12){\line(1, 0){160}}
\put (-65, 0) {$ \exists \, K \mr{f} C \; \; such \; that \;  \;  \sigma = uf$}
 \end{picture}
 \end{center}

\vspace{1ex}
 
\item \label{y3}   
Given any quasispace  $(S, \, X)$, there is an equivalence
\mbox{(natural in $K$)} 
 \begin{center}
\begin{picture}(0, 30) 
\put (-50, 18)  {$uK \mr{\sigma} S \; \in \;  XK$}
\put (-135, 12){\line(1, 0){260}}
\put (-120, 0) {$(uK,\, \varepsilon K) \mr{\sigma} (S, \, X) \; is \; a
  \; morphism \; of \; quasispaces$}
 \end{picture}
 \end{center}

\vspace{1ex}

\noindent That is, there is an equality of sets $XK = [\varepsilon K,
  \, (S, \, X)]$.  

\vspace{1ex}

\item \label{y4}

\noindent If covers are final families in $\bb{C}$, then $\varepsilon$
is full.

\noindent If $u$ is faithful, then  $\varepsilon$ is faithful.

\vspace{1ex}

\item \label{y5}

The family  $\varepsilon C  \mr{\sigma}  (S, \, X)$, $\sigma
\in XC$,  is a final family in $\bb{Q}$.

\vspace{1ex}

\item \label{y6}

Given any cover $K_i \to K \in \cc{J}_K$, the family 
$\varepsilon K_i \to \varepsilon K$ is a final family in $\bb{Q}$. If
the cover is surjective, it is a final surjective (thus strict
epimorphic) family in  $\bb{Q}$. 
  
\newcounter{yoneda} 
\setcounter{yoneda}{\theenumi}
\end{enumerate}   
 
 \vspace{1ex}

{\bf proof:} To check items (\ref{y1}) and (\ref{y3}) is sharp but straightforward, and
 it is left to the reader. Item (\ref{y2}) is immediate from item
 (\ref{y1}), and  item
 (\ref{y5}) from item  (\ref{y3}). Item (\ref{y4}) follows easily from
 items (\ref{y2}) and (\ref{y3}). Finally, to check item (\ref{y6}) is
 again sharp but
 straightforward, or, if the reader prefers, it follows immediately
 from \mbox{proposition \ref{charfinalinQ}} below. 
 \qed\end{yoneda}

\vspace{1ex}
 
Initial families in $\bb{Q}$ are easily characterized. The proof of
the following proposition is immediate:

\vspace*{4ex}

\begin{proposition} \label{charinitialQ}
 A family $(S,\, X) \mr{\varphi_\alpha} (S_\alpha , \, X_\alpha)$ in $\bb{Q}$ is initial if and only if given any $C$ and $uC \mr{\sigma} S$, the following equivalence holds:
\begin{center}
\begin{picture}(0, 35) 
\put (-38, 18)  {$uC \mr{\sigma} S \; \in\; XC$}
\put (-80, 12){\line(1, 0){160}}
\put (-65, 0) {$\forall \alpha \;\;\;  uC \mr{\sigma} S \mr{\varphi_\alpha} S_\alpha  \; \in\; X_\alpha C$}
 \end{picture}
 \end{center}
\qed \end{proposition}
\begin{theorem} \label{Qistopological}
The functor $\bb{Q} \mr{q} \bb{S}$ is a topological functor. 
\end{theorem}
\begin{proof}
We have to see that $q$ creates initial families (see
\ref{topological} (\ref{topin})). Given quasispaces  $(S_\alpha , \, X_\alpha)$ and a family $S \mr{\varphi_\alpha} S_\alpha$ in $\bb{S}$, define a quasispace structure $X$ on $S$ stipulating that $uC \mr{\sigma} S \; \in\; XC$ by the equivalence in proposition \ref{charinitialQ}. It is immediate to check that the pair $(S, \, X)$ so defined is a quasispace and that the $\varphi_\alpha$ become morphisms of quasispaces $(S,\, X) \mr{} (S_\alpha , \, X_\alpha)$. 
\end{proof}

It follows from \ref{charinitialQ} and \ref{topological}
(\ref{topselfdual}) that $u$ creates final families. However, we shall
prove this directly because the proof yields an essential
characterization of final families reminiscent of the characterization
of epimorphic families of sheaves. 

\begin{proposition} \label{charfinalinQ}
A family $(S_\alpha , \, X_\alpha) \mr{\varphi_\alpha} (S,\, X)$ in $\bb{Q}$ is final if and only if given any $C$ and $uC \mr{\sigma} S$, the following equivalence holds:
$$
\xymatrix@R=5pt
               {
                \\
                \sigma \in XC \;\iff\; \exists \, C_i \to C \in \cc{J}_C, \; 
                \sigma_i \in X_{\alpha_i}\,,  \; such  \; that \;  
                }
\xymatrix@R=15pt@C=20pt
              {
                   uC_i \ar[r]^{\sigma_i}  \ar[d]  
               &  S_{\alpha_i}  \ar[d]^{\varphi_{\alpha_i}}
               \\
                  uC \ar[r]^{\sigma}  
               & S
              }
$$
\end{proposition} 
\begin{proof}
Given quasispaces  $(S_\alpha , \, X_\alpha)$ and a family $S_\alpha \mr{\varphi_\alpha} S$ in $\bb{S}$, define a quasispace structure $X$ on $S$ stipulating that $uC \mr{\sigma} S \in XC$ by the equivalence above. To prove the statement we have to check that $(S, \, X)$ is a quasispace, that the $\varphi_\alpha$ become morphisms of quasispaces $(S_\alpha,\, X_\alpha) \mr{} (S, \, X)$, and that the resulting family is final. 

1) $(X,\, S)$ is a quasispace:

\noindent \emph{presheaf condition}:
                  Given  $C \mr{f} D$ and  $uD \mr{\sigma} S \, \in XD$, 
                 take  $C_j  \to C  \, \in \cc{T}_C$
                 and  
$$
\xymatrix@R=0.1pt
                {
                 \\                 
                 C_j \mr{f_j} D \;\;  such \;  that 
                }
\xymatrix@R=15pt@C=20pt
              {
                   C_j \ar[r]^{f_j}  \ar[d]  
               &  D_{i_j}  \ar[d]
               \\
                  C \ar[r]^{f}  
               & D 
              }
\xymatrix@R=5pt
               {
                \\
                 . \; We \; have \; then
               }
\xymatrix@R=15pt@C=20pt
              {
                   uC_j       \ar[r]^{uf_j}             \ar[d]  
               &  uD_{i_j} \ar[r]^{\sigma_{i_j}}  \ar[d]  
               &  S_{\alpha_{i_j}}                      \ar[d]^{\varphi_{\alpha_{i_j}}}
               \\
                  uC   \ar[r]^{uf} 
               & uD   \ar[r]^{\sigma}
               & S  
              }
$$
(we have used here property  (U) of $\cc{T}$). From the presheaf condition on $X_{\alpha_{i_j}}$ it follows that $\sigma_{i_j} \circ uf_j \,\in\, X_{\alpha_{i_j}} C_j$. This shows that $\sigma \circ uf \,\in\, XC$.
 
 \vspace{1ex}

\noindent \emph{covering condition}: Let $C_i \to C  \, \in \cc{T}_C$ and 
$uC \mr{\sigma} S$  be such that the composites $(uC_i \to uC \mr{\sigma} S) \in XC_i$. Take (for each $i$)  a cover $C_{i, j} \to C_i \, \in \cc{T}_{C_i}$ and maps 
$$
\xymatrix@R=3pt
               {
                \\
               uC_{i, j} \mr{\sigma_{i, j}}  S_{\alpha_{(i, j)}} \; \in \, X_{\alpha_{(i, j)}} C_{i, j} 
                \; \; such \; that \;
               }
\xymatrix@R=15pt@C=20pt
               {
                       uC_{i, j} \ar[rr]^{\sigma_{i, j}}  \ar[d]  
                &&  S_{\alpha_{(i, j)}} \ar[d]^{\varphi_{\alpha_{(i, j)}}}
                \\
                       uC_i  \ar[r] 
                &     uC    \ar[r]^{\sigma}
                &    S   
               } 
$$
This shows that $\sigma \in XC$ (we use now property (C) of $\cc{T}$).

\vspace{1ex}

2)  $(S_\alpha,\, X_\alpha) \mr{\varphi_\alpha} (S, \, X)$ is a quasispace morphism:
Given $\alpha$ and a map 
$
\xymatrix@R=3pt
               {
                \\
                 uC \mr{\sigma} S_\alpha \, \in X_\alpha C, \; consider \; the \; diagram:\; 
               }
\xymatrix@R=15pt@C=30pt
               {
                     uC \ar[r]^\sigma  \ar[d]^{id}
                &   S_\alpha \ar[d]^{\varphi_\alpha}
                \\
                     uC \ar[r]^{\varphi_\alpha \circ \sigma} 
                &  S     
               } 
\xymatrix@R=3pt
               {
                \\
                . \; This \; shows \; that
                }
$               
  \mbox{$\varphi_\alpha \circ \sigma  \in XC$} 
(we use now property (I) of $\cc{T}$).  

\vspace{1ex}

3) Finally, to check that the resulting family is a final family is straightforward.         
\end{proof}
 An important consequence of this characterization is that final families in $\bb{Q}$ are universal. We have:
 \begin{proposition} \label{qfinalareuniversal}
 Final families for the functor $\bb{Q} \mr{q} \bb{S}$ are universal (property  (U) in definition \ref{familyproperties}).
 \end{proposition}
 \begin{proof}
 Let $(S_\alpha , \, X_\alpha) \mr{\varphi_\alpha} (S,\, X)$ be a final family and 
 $(T, \, Y) \mr{\phi} (S, \,X)$ a morphism of quasispaces. Consider the  $\bb{C}$-crible on $T$ defined by: 
 $$uC \mr{\varphi} T  \in  P  \iff $$
 $$
 \xymatrix@R=0pt
              {
               \\
                \varphi \in YT, \; and  \; \;\exists \alpha \;\; and\;\;
                 uC \mr{\theta} S_\alpha , \; \theta \in X_\alpha C, \; \; \;such\;that                               }               
\xymatrix@R=15pt@C=20pt
        {
         uC \ar[r]^\theta \ar[d]^{\varphi} 
        & S_{\alpha} \ar[d]^{\varphi_\alpha}
        \\
        T \ar[r]^\phi 
        & S
        }
$$

This defines an $r$-pullback in $\bb{Q}$ by furnishing $uC$ with the
 yoneda $\varepsilon C$ quasispace structure (\ref{yoneda}). We shall show that the family
 $(uC, \, \varepsilon C) \mr{\varphi}  (T, \, Y)$, $\varphi \in P$, is a final family in $\bb{Q}$. 
 
 Given $uK \mr{\eta} T \in \, YK$, the composite $\phi \circ \eta \in \, XK$. Using  proposition \ref{charfinalinQ} take  $K_i \mr{s_i} K \in \  \cc{J}_K$ and 
 $uK_i \mr{\eta_i} S_{\alpha_i} \in \, X_{\alpha_i}K_i$ such that in the following diagram the exterior commutes 
 $$
 \xymatrix@R=15pt@C=20pt
                 {
                      uK_i   \ar @<+2pt> `u[r] `[rr]^{\eta_i} [rr]
                                \ar[r]^{id}  \ar[d]^{us_i}
                 &   uK_i  \ar[d]^{\varphi_i} \ar[r]^{\eta_i}
                 &   \;\; S_{\alpha_i}  \ar[d]^{\varphi_{\alpha_i}}
                 \\
                      uK  \ar @<+2pt> `d[r] `[rr]_{\phi \circ \eta} [rr]
                            \ar[r]^\eta
                 &   T \ar[r]^\phi
                 &   S
                 } 
$$
Fill in the middle vertical arrow, with $\varphi_i = \eta \circ us_i$. The square on the right shows that $uK_i \mr{\varphi_i} T \, \in \, P$. Then, the square on the left finishes the proof by another application of proposition \ref{charfinalinQ} .
\end{proof}

Actually we are interested in final surjective families. Concerning surjective families we have:

\begin{proposition}
If strict epimorphic families in $\bb{S}$ are universal, then surjective families for the functor $\bb{Q} \mr{q} \bb{S}$ are universal (property  (U) in definition \ref{familyproperties}). 
\end{proposition}
\begin{proof}
Just observe that the initial structure determined by the two arrows out of the upper left corner  of an $r$-pullback (definition \ref{rpullback}) taken in $\bb{S}$ yields a 
$r$-pull-back in $\bb{Q}$.
\end{proof}

Since surjective and final families have property (S), it follows from
 \ref{collectionfacts} (\ref{intersection}) that when strict
 epimorphic families in $\bb{S}$ are universal, final surjective
 families in 
 $\bb{Q}$ are universal. Thus, from theorem \ref{Qistopological} and
 corollary \ref{tisq2} we have: 

\begin{theorem} \label{qisquasi}
If strict epimorphic families in $\bb{S}$ are universal (in
particular, if $\bb{S}$ is $f$-regular),  then so they are in
$\bb{Q}$, and the functor $\bb{Q} \mr{q} \bb{S}$ is $f$-regular.
\end{theorem}
 
 Then, from proposition \ref{quasi=fmande}  and theorem \ref{lift} we have:

\begin{theorem} [compare \cite{P2}, 5.8] \label{Qisaquasitopos}
A category $\bb{Q} \mr{q} \bb{S}$ of quasispaces over a $f$-regular category  $\bb{S}$ is $f$-regular, and if $\bb{S}$ is a quasitopos, then so it is $\bb{Q}$.
\end{theorem}

\vspace{1ex}

Let $\bb{P}$ be the category of quasispaces for the  trivial (generated by the
isomorphisms in $\bb{C}$) Grothendieck topology. Clearly, the
inclusion determines
full and faithful functor $\bb{Q} \mmr{c} \bb{P}$, where $\bb{Q}$ is
the category of quasispaces determined by any other topology on
$\bb{C}$. We shall construct a left adjoint $\bb{P}
\mr{\#} \bb{Q}$ to the
functor $c$, and study its basic properties. We consider the more
general situation determined by an inclusion of Grothendieck
topologies.

\vspace{1ex}

\begin{sinnada} \label{aqf} {\bf Associate Quasispace Functor.} Consider
a category  $\mathbb{C}$, 
two Grothendieck topologies  $\cc{T} \subset \cc{J}$, and 
a functor \mbox{$\mathbb{C} \mr{u} \mathbb{S}$.} Let $\bb{P} \mr{p}
\bb{S}$, $\bb{Q} \mr{q} \bb{S}$
be the respective categories of quasispaces, and $\bb{C} \mr{h} \bb{P}$, 
$\bb{C} \mr{\varepsilon} \bb{Q}$ the respective yoneda functors
(\ref{yoneda}). Let  $\bb{Q} \mmr{c} \bb{P}$ be the inclusion
functor (notice that $p \circ c = q$, but in general $c \circ
\varepsilon \neq h$).
\end{sinnada} 
 We have:
\begin{proposition}
The inclusion functor $\bb{Q} \mmr{c} \bb{P}$ has a left adjoint,
$\# \dashv c$, \mbox{$\bb{P} \mr{\#} \bb{Q}$,} $id \mono
c \circ \#$, $\# \circ c = id$. Given $(S, \, X) \in \bb{P}$, $\#(S, \,
X) = (S, \, \#X)$ is defined by stipulating the equivalence:
\begin{center}
\begin{picture}(0, 30) 
\put (-50, 18)  {$uK \mr{\sigma} S \; \in \; \#X(K)$}
\put (-115, 12){\line(1, 0){220}}
\put (-110, 0) {$\exists \; K_\alpha \mr{} K \in \, \cc{J}_K \;|\;\;
  uK_\alpha \mr{} uK \mr{\sigma} S \;\in \;XK$} 
 \end{picture}
 \end{center}
Moreover,   $q \circ \# = p$ and $\# \circ h = \varepsilon$.
\end{proposition}
\begin{proof}
The proof is sharp but straightforward and it is left to the
reader. We mention  that the properties (I), (C) and
(U) (see definition \ref{familyproperties}) of the covers are essential.
\end{proof}

\begin{proposition}
The functor $\bb{P} \mr{\#} \bb{Q}$ preserves finite strict
monomorphic \mbox{families.} 
\end{proposition}
\begin{proof}
Since $p$ and $q$ are, in particular, $\cc{M}_{fin}$-functors, it is
equivalent to work with initial injective families. We shall see that,
in fact, $\#$ preserves any (not necessarily injective) finite
initial family. Clearly $\#$ preserves the empty initial family
$S_\top$ since this 
family is a quasispace for any topology (example \ref{example1}). Let now
$(S,\,X) \to (S_i , \, X_i)$ be a finite initial family. The following
chain of equivalences proves the proposition:
\begin{center}
\begin{picture}(0, 100) 
\put (-75, 88)  {$\forall \, i \;\; uK \mr{\sigma} S \mr{} S_i \; \in \; \#X_i(K)$}
\put (-150, 80){\line(1, 0){300}}
\put (-145, 66) {$\forall \, i \; \exists \; K_{i,\alpha} \mr{} K
  \in \, \cc{J}_K  \;\; | \;\;   uK_{i,\alpha} \mr{} uK \mr{\sigma} S
  \mr{} S_i \;\in \;X_i(K_{i,\alpha})$}
\put (-160, 56){$a$} \put (-150, 58){\line(1, 0){300}}
\put (-138, 44)  {$ \exists \, K_\lambda \mr{} K  \in \, \cc{J}_K \;\; |\;\; \forall \,
  i \;\; uK_\lambda \mr{} uK \mr{\sigma} S
  \mr{} S_i \;\in \;X_i(K_{\lambda})$}
\put (-150, 36){\line(1, 0){300}}
\put (-115, 22)  {$\exists \, K_\lambda \mr{} K \,  \in \, \cc{J}_K
  \;\; | \;\; uK_\lambda
  \mr{}  uK \mr{\sigma} S \; \in \; X(K_{\lambda})$}
\put (-150, 14){\line(1, 0){300}}
\put (-50, 0)  {$uK \mr{\sigma} S \; \in \; \#X(K)$}
\end{picture}
 \end{center}
All unlabeled equivalences hold by definition, and the equivalence $a$
holds by \mbox{property (F)} (filtered, definition \ref{familyproperties}) of
the covers.
\end{proof}

\begin{theorem} \label{assisfregular}
If strict epimorphic families in $\bb{S}$ are universal (in particular, if $\bb{S}$ is $f$-regular),  then the functor $\bb{P} \mr{\#} \bb{Q}$ is $f$-regular.
\end{theorem}
\begin{proof}
Recall that $q \circ \# = p$. We refer to definition
\ref{fregularfunctor2}. Conditions  RF1) and RF2) are clear. By RF1)
we know that $\#$ preserves $s\cc{E}$-families, and in the previous
proposition we established that it preserves
$s\cc{M}_{fin}$-families. By the assumption made,  $p$ and $q$ are $f$-regular
functors (theorem \ref{qisquasi}), so conditions RF3) and RF4) follow
from proposition \ref{efunctoryoga1} (and RF5) is equivalent to RF4)). 
\end{proof}

\vspace{2ex}

\section{Strict quasispaces over a $f$-regular category $\bb{S}$} \vspace{1ex}
 In the classical examples  a third condition is required in the
 definition of quasispace. In this general setting it means a
 \emph{density} condition for the admissible maps. Nothing much can be
 proved if we do not put some restrictions on the category $\bb{S}$.
Although not always necessary, the sensible generality is to assume
 that $\bb{S}$ is $f$-regular.
\begin{assumption} \label{data2}
 The category $\bb{S}$ is a  $f$-regular category \emph{(in particular strict epimorphic families in $\bb{S}$ have all five properties in definition  \ref{familyproperties}).}  
\end{assumption}
 
\begin{definition} \label{strictquasispace}
A \emph{strict quasispace} is a quasispace which in addition satisfies the following condition:

Density condition:
$$
The \; family \;\; UC \mr{\sigma} S , \; \sigma \in XC,\; all \; C \in \bb{C}, \;\; is \; strict \; epimorphic  \; in \; \bb{S}. 
$$
\end{definition}
The category $s\bb{Q}$ is defined as the full subcategory of the
category $\bb{Q}$ of quasispaces whose objects are the strict
quasispaces. We denote by $q_s$ the restriction of the functor $q$. 
There is a commutative diagram of faithful functors: 
$$
\xymatrix@R=15pt@C=7pt
         {
          s\mathbb{Q}\; \ar@{^{(}->}[rr]^{i}  \ar[rd]_{q_s}  
           &&  \mathbb{Q} \ar[ld]^{q}  
           \\
           & \bb{S} 
         }
$$ 
\begin{yoneda}[continued] $\,$ \label{yonedacontinued}
\begin{enumerate} 
\setcounter{enumi}{\theyoneda}

\item \label{y7} 
For any object $C \in \bb{C}$, the quasispace $\varepsilon C$ is
strict. Thus there is a factorization (that we also denote
$\varepsilon$)
$$        
\xymatrix@1
        {
         \bb{C} \;\; \ar @<+2pt> `u[r] `[rr]^{\varepsilon} [rr]
                                             \ar[r]^{\varepsilon} 
         & \;\;s\bb{Q}\;\;  \ar@{^{(}->}[r]^{i} 
         & \;\;\bb{Q}      
        }
$$

\item \label{y8}
Given any strict quasispace $(S, \, X)$, the family  $\varepsilon C
\mr{\sigma}  (S, \, X)$, \mbox{$\sigma \in XC$,}  is a final
surjective (thus strict epimorphic) family in $\bb{Q}$ (thus also in
$s\bb{Q}$)
\setcounter{yoneda}{\theenumi}
\end{enumerate}

{\bf proof:} Item (\ref{y7}) is clear since $id_{uC}  \in  \varepsilon
C$). Item (\ref{y8}) holds by definition of strict quasispace and item
(\ref{y5}) in \ref{yoneda}. 
\qed \end{yoneda}

Any $f$-regular functor $\bb{C} \mr{u} \bb{S}$ into a $f$-regular
category $\bb{S}$ is the forgetful functor of a category of strict
quasispaces. We have:

\begin{theorem}[compare \cite{P2}, 5.13]
Given  a $f$-regular functor $\bb{C} \mr{u} \bb{S}$, with $\bb{S}$
$f$-regular. Consider the Grothendieck topology in $\bb{C}$ whose
covers are the strict epimorphic  families. Then, in the commutative
triangle
$
\xymatrix@R=15pt@C=7pt
         {
          \bb{C} \; \ar[rr]^{\varepsilon}  \ar[rd]_{u}  
           &&  s\mathbb{Q} \ar[ld]^{q_s}  
           \\
           & \bb{S} 
         }
$, 
the functor $\varepsilon$ is an equivalence.
\end{theorem}
\begin{proof}
From  proposition \ref{quasi=fmande} we know that $\bb{C}$ is
$f$-regular, so that the strict epimorphic families are a
Grothendieck topology. By yoneda \ref{yoneda} (\ref{y4}) 
$\varepsilon$ is full and faithful. Let now $(S, \, X)$ be any strict
quasispace, and consider a strict epimorphic family $K \mr{f_\sigma}
C$ in $\bb{C}$ over the strict epimorphic family $uK \mr{\sigma} S$,
all $K \in \bb{C}$ and $\sigma \in XK$. Since $f_\sigma$ is a cover,
by yoneda \ref{yoneda} (\ref{y6}) the family $\varepsilon K
\mr{\varepsilon f_\sigma} \varepsilon C$ is strict epimorphic, and
by  yoneda \ref{yonedacontinued} (\ref{y8}) the family  
$\varepsilon K \mr{\sigma} (S,\,X)$ is also strict epimorphic and sits
over  $uK \mr{\sigma} S$. It follows that $\varepsilon C$ is
isomorphic to $(S,\,X)$ over $id_S$.
\end{proof}
\vspace{1ex}

We see in particular that in this case the functor $q_s$ is
$f$-regular. We study now the $f$-regularity condition for the functor $q_s$ in
general.  

Since strict epimorphic families compose in $\bb{S}$, it easily follows:
\begin{proposition} \label{strictsurjective}
Given a surjective family of quasispaces $(S_\alpha ,\, X_\alpha) \mr{\varphi_\alpha} (S, \, X)$, if each  $(S_\alpha ,\, X_\alpha)$ is strict, then so it is  $(S, \, X)$.
\qed\end{proposition}

\begin{proposition} \label{strictepiofstrict}
Given a strict epimorphic family of quasispaces \mbox{$(S_\alpha ,\, X_\alpha)
\mr{\varphi_\alpha} (S, \, X) \;\in\; \bb{Q}$,} if each  $(S_\alpha ,\, X_\alpha)$ is
strict, then so it is  $(S, \, X)$. 
\end{proposition}
\begin{proof}
Since $\bb{Q} \mr{q} \bb{S}$ is topological (thus an
$\cc{E}$-functor), the given family is final surjective
(\ref{characterization2}). Then,  by the previous proposition $(S, \,
X)$ is strict. 
\end{proof}

Thus we have the following:

\begin{corollary} \label{sandiareEfunctors}
 $s\bb{Q}$ is closed under strict epimorphic families in $\bb{Q}$, and
  the functors $q_s$ and $i$ are $\cc{E}$-functors 
\qed\end{corollary}

\begin{proposition} \label{qshasright}
The functor $s\bb{Q} \mr{q_s} \bb{S}$ has a right adjoint 
$q_s \dashv r$, 
 $\bb{S} \mr{r} s\bb{Q}$, and if for any  $S \in \bb{S}$ the family  $uC \mr{\sigma} S$, all $C \in \bb{C}$, is
 strict  epimorphic, then the right adjoint is full and
faithful, $r = (-)_\top$, $q_s (-)_\top = id$.
\qed\end{proposition} 
\begin{proof} Follows by proposition \ref{preservesradjoint1}. Given $S \in \bb{S}$, take the
 family $uC \mr{\sigma} S$, all $C \in \bb{C}$,
factor this family 
\mbox{$uC \mr{\psi} H \mr{m}  S$,} with a monomorphism 
$m$ and the family $\psi$ strict epimorphic. Let
 $rS = (H, rS)$ be a strict epimorphic
  family over the family
$\psi$.  Then, $rS$ is the right adjoint to $q_s$ on $S$. The second
 claim is clear since in this case $H = S$.
\end{proof}

Notice that the condition in the proposition says that
the quasispace $S_\top$ is strict. When this in not the case, there
can not be any strict quasispace structure on $S$, that is, the fiber
$s\bb{Q}_S$ is empty.

\begin{proposition} \label{ihasarightadjoint}
The inclusion $s\bb{Q} \mr{i} \bb{Q}$ has a right adjoint  $\bb{Q} \mr{s} s\bb{Q}$,  
$i \dashv s$ (notice that  $ r  = s(-)_\top$ and $\; q_s s \neq q$).
\end{proposition}
\begin{proof}
Follows by proposition \ref{preservesradjoint1}. Given a quasispace $(S,\, X)$, we abuse notation and
denote $s(S, \, X) = (sS, \, X)$, where $sS \subset S$ is given by the
factorization of the family of all admissible maps $uC \mr{\sigma} S$, $C \in \bb{C}$ into a strict epimorphic family followed by a monomorphism 
\mbox{$uC \mr{\sigma} sS \subset S$.}
\end{proof}

We see from proposition \ref{strictsurjective} that a surjective family of strict quasispaces $(S_\alpha ,\, X_\alpha) \mr{\varphi_\alpha} (S, \, X)$ is final in $s\bb{Q}$ if and only if it is characterized by the equivalence in proposition
\ref{charfinalinQ}. This is not the case for non surjective final
families in $s\bb{Q}$. 
On spite of having a right adjoint, the inclusion functor $i$ will not preserve final families unless they are surjective. It is clear that the quasispace $S_\bot$ (the empty final family, see example \ref{example1}) is not a strict quasispace. More generally, at least when the covers are surjective, we have:

\begin{remark} \label{finalonlyifsurjective}
Any family $(S_\alpha ,\, X_\alpha) \mr{\varphi_\alpha} (S, \, X)$ between strict quasispaces which is final in $\bb{Q}$ (thus characterized by the equivalence in proposition \ref{charfinalinQ}) is necessarily surjective (assuming the covers to be surjective).  
\end{remark}
\begin{proof}
For each $UC \mr{\sigma} S$,  $\sigma \in XC$, $C \in \bb{C}$, take
a cover $C_{\sigma, i} \to C \, \in \cc{T}_{C}$ and maps 
$$
\xymatrix@R=3pt
               {
                \\
               uC_{\sigma, i} \mr{}  S_{\alpha_{(\sigma, i)}} \; \in \, X_{\alpha_{(\sigma, i)}} C_{\sigma, i} 
                \; \; \; such \; that \;\;
               }
\xymatrix@R=15pt@C=20pt
               {
                    uC_{\sigma, i} \ar[r]  \ar[d]  
                &  S_{\alpha_{(\sigma, i)}} \ar[d]^{\varphi_{\alpha_{(\sigma, i)}}}
                \\
                       uC  \ar[r]^\sigma
                &     S                 } 
$$
By property (C) the composite by the lower left corner is strict epimorphic. Then the statement follows by property (S).  
\end{proof}

\vspace*{7ex}

\begin{proposition} \label{q_sfinalareuniversal}
If covers are surjective families in $\bb{C}$, final surjective families for the functor $s\bb{Q} \mr{q_s} \bb{S}$ are universal (property  (U) in definition \ref{familyproperties}).
 \end{proposition}
\begin{proof}
As remarked above, by proposition \ref{strictsurjective} a final
surjective family in $s\bb{Q}$ is characterized by proposition
\ref{charfinalinQ}. Then it follows from \ref{yonedacontinued}
(\ref{y7}) that the r-pullback constructed in the proof of proposition \ref{qfinalareuniversal} is an r-pullback of strict quasispaces. The resulting family is a final family in $\bb{Q}$, then by remark 
\ref{finalonlyifsurjective} it is also surjective.
\end{proof}

\emph{The fact that final surjective families in $s\bb{Q}$ are universal
holds independently of whether the covers are surjective families in
$\bb{C}$, or not.}

However, some other assumptions have to be made. These
assumptions are forced upon us by the need that a \emph{finite initial
  injective} family of strict quasispaces  be in fact an strict
quasispace. Given strict quasispaces  $(S_\alpha , \, X_\alpha)$ and a
family $S \mr{\varphi_\alpha} S_\alpha$ in $\bb{S}$, the initial
quasispace structure defined on $S$ in theorem \ref{Qistopological}
will not be strict, unless we make further assumptions. The failure or
not of this structure to be strict is related to the fact of whether
the functor $q_s$ is $f$-regular (that implies, in particular, that
final surjective families in $s\bb{Q}$ are universal). When \emph{any
  initial injective} family of strict quasispaces is a strict
quasispace, the functor $q_s$ will be topological.

\begin{assumption} \label{a2}
Given a finite initial injective family of quasispaces 
\mbox{$(S, \, X) \mr{\varphi_\alpha} (S_\alpha ,\, X_\alpha)$,} if each  $(S_\alpha ,\, X_\alpha)$ is strict, then so it is  $(S, \, X)$.
\end{assumption}

\begin{proposition}[Under assumption \ref{a2}]
Given a finite strict monomorphic family of quasispaces
\mbox{$(S, \, X) \mr{\varphi_\alpha} (S_\alpha ,\, X_\alpha) \;\in\;
  \bb{Q}$,} if each  $(S_\alpha ,\, X_\alpha)$ is 
strict, then so it is  $(S, \, X)$. 
\end{proposition}
\begin{proof}
Since $\bb{Q} \mr{q} \bb{S}$ is topological (thus a
$\cc{M}_{fin}$-functor), the given family is initial injective
(\ref{finitecharacterization2}). Then,  by the assumption $(S, \,
X)$ is strict. 
\end{proof}

Thus we have the following:

\begin{corollary} \label{sandiareMfunctors}
 $s\bb{Q}$ is closed under finite strict monomorphic families in $\bb{Q}$, and
  the functors $q_s$ and $i$ are $\cc{M}_{fin}$-functors 
\qed\end{corollary}

 In particular, $s\bb{Q}$ is closed under pullbacks taken in $\bb{Q}$
 (\mbox{proposition \ref{preservescolandlim}).}
\begin{theorem}[Under assumption \ref{a2}] \label{q_sisquasitopological}
The functors 
$s\bb{Q} \mr{q_s} \bb{S}$ and $s\bb{Q} \mr{i} \bb{Q}$ are $f$-regular.
\end{theorem}
\begin{proof}
We already know that they have a right adjoint (propositions
\ref{qshasright} and \ref{ihasarightadjoint}) and that they are
$\cc{E}$ and $\cc{M}_{fin}$-functors (propositions
\ref{sandiareEfunctors} and \ref{sandiareMfunctors}). To prove that they are $f$-regular it remains to see that strict epimorphic families in $s\bb{Q}$ are universal. But $\bb{Q}$ is a $f$-regular category (theorem 
\ref{Qisaquasitopos}), and the inclusion $s\bb{Q} \mr{i} \bb{Q}$ is
closed under pullbacks and strict epimorphic families. The claim follows. 
\end{proof}

Then, from proposition \ref{quasi=fmande}  and theorem \ref{lift} we have:

\begin{theorem}[Under assumption \ref{a2}), (compare \cite{P2}, 5.10]
  \label{sQisaquasitopos} 
A category $s\bb{Q} \mr{q} \bb{S}$ of strict quasispaces over a $f$-regular category  $\bb{S}$ is $f$-regular, and if $\bb{S}$ is a quasitopos, then so it is $s\bb{Q}$.
\qed\end{theorem}


From theorems \ref{assisfregular} and \ref{q_sisquasitopological} it follows

\begin{sinnada} {\bf Associate Quasispace Functor} \label{aqfcontinued}
  (continued). Consider the situation described in \ref{aqf}. It is
  clear that if $(S,, X) \in \bb{P}$ is an strict quasispace, then so
  it is  $\#(S,, X) = (S,, \#X) \in \bb{Q}$, and we have a functor
 $s\bb{P} \mr{\#} s\bb{Q}$ left adjoint to the inclusion $s\bb{Q} \mr{c}
  s\bb{P}$. There is the following square of pairs of adjoint
  functors:
$$
\xymatrix@R=10pt
         {{} \\ {\# \dashv c\, , \;\; i \dashv s\, ,}}
\hspace{10ex} 
\xymatrix@R=25pt@C=30pt
              {
                   s\bb{Q} \ar@<-1ex>[r]^{c}  \ar@<-1ex>[d]_{i}  
               &  s\bb{P}  \ar@<-1ex>[l]_{\#} \ar@<-1ex>[d]_{i}
               \\
                  \bb{Q} \ar@<-1ex>[u]^{s} \ar@<-1ex>[r]^{c}   
               & \bb{P} \ar@<-1ex>[u]^{s} \ar@<-1ex>[l]_{\#}
              }
$$ 
The inclusion functors commute, and the
adjunctions in the full subcategories are the restrictions of the
adjunctions in the larger categories 
$$\; i \circ c = c \circ i\,,  \;\;\;   i \circ \# = \# \circ i\,,
\;\;\;  s \circ c = c \circ s\,.
$$
\emph{(Under assumption \ref{a2})}. The functors $i$, $\#$, are
$f$-regular functors (the 
  $f$-regularity for the functor $\#$ between strict quasispaces
  follows because the inclusion functors $i$ are $f$-regular). 
\qed\end{sinnada}

The following conditions were considered by Penon (\cite{P2}, 5.10)  to
prove that the category $s\bb{Q}$ is a quasitopos. We use them to
insure the validity of assumption \ref{a2}.

\begin{proposition} \label{sufficient}
Let $\bb{C} \mr{u} \bb{S}$ $(\,\bb{S}$ $f$-regular $)$  be such that:

i) The quasispace $1_\top$ is strict.

ii) Given any $C \in \bb{C}$ and a strict subobject $S \mono uC =
q \varepsilon C$, the initial quasispace structure induced on $S$ is
strict.

iii) Given any two objects $C, \; D \in \bb{C}$, the product 
$\varepsilon C \times \varepsilon D$ taken in  $\bb{Q}$ is an strict quasispace.

Then, finite initial injective families of strict quasispaces are
strict, that is, the assumption \ref{a2} holds.
\end{proposition}
\begin{proof}
The empty initial injective family is given by $1_\top$. We claim that given a
strict quasispace $(S,\,X)$ and a strict monomorphism $T \mono S$, the
initial quasispace $(T,\,Y) \mono (S,\,X)$ is strict: We argue over
the following diagram:
$$ 
\xymatrix@R=20pt@C=20pt
          {
             \varepsilon K \ar[r]^\mu  \ar[rd]^{\!\!\!\!\sigma_{\mu, \tau}}
           & (P_\tau,\, Z_\tau)  \ar@{^{(}->}[r]  \ar[d]^{\pi_\tau}
           & \varepsilon C \ar[d]^\tau
           \\
           & (T,\,Y) \ar@{^{(}->}[r] 
           & (S,\,X)
          }
$$                      
The family $\{\tau, \; all \; C \in \bb{C},\; \tau \in XC\}$ is
strict epimorphic (yoneda \ref{yonedacontinued} (\ref{y8})). Pulling
back this family in $\bb{Q}$ we have the 
family $\pi_\tau$, which is strict epimorphic since $\bb{Q}$ is
$f$-regular (theorem \ref{Qisaquasitopos}). Since $Z_\tau$ is the
initial structure induced by $P_\tau \mono uC$ (the reader can
check this), for each $\tau$, by ii),
the family $\{\mu, \; all K \in \bb{C}, \mu \in Z_{\tau}K\}$ is strict
epimorphic. The composite family $\sigma_{\mu, \tau} = \pi_\tau \circ
\mu$ is strict epimorphic. But $\sigma_{\mu, \tau} \in YK$, the claim
follows. 

We let the reader verify that from iii) it easily follows that
the product of two (hence any finite product) of strict quasispaces is
strict.   To finish the proof recall that a initial injective
family induces a initial injective map into a product. 
\end{proof}
  
\vspace{1ex}

In general the functor $s\bb{Q} \mr{q_s} \bb{S}$ will not be topological.  It is clear that the quasispace $S_\bot$ is not strict, but still it may exist a smallest strict quasispace structure on $S$. This is the case in many classical examples (see \cite{D}), and it is equivalent to the fact that the functor $s\bb{Q} \mr{q_s} \bb{S}$ is topological. We consider a general situation that include all these examples:

\begin{assumption} \label{classical} 
There is a class of objects $\bb{I} \subset \bb{C}$ such that:

a) For each $I$ in $\bb{I}$ and $C$ in $\bb{C}$, $u$ establish a bijection
\begin{picture}(0, 18) 
\put (12, 5)  {$I \mr{} C \; \; in\; \bb{C}$}
\put (2, 2){\line(1, 0){82} \;.}
\put (7, -7) {$uI \mr{} uC \;\;  in\;\bb{S}$}
 \end{picture}

\vspace{2ex}
 
b) For any  $S$ in $\bb{S}$, the family of all $\; uI  \to S$, $I \in \bb{I}$, is strict epimorphic.

\vspace{1ex}

c) Given any strict epimorphic family $S_\alpha  \to S$ in $\bb{S}$, any $I \in \bb{I}$ and $uI \mr{\sigma} S$,

\vspace{1ex}

$
\xymatrix@R=5pt
               {
                \\
                there \; exists \; I_i \to I \in \cc{T}_I, \; 
                u I_i  \to  S_{\alpha_i},  \; such  \; that \;  
                }
\xymatrix@R=15pt@C=20pt
              {
                   uI_i \ar[r]  \ar[d]  
               &  S_{\alpha_i}  \ar[d]^{\varphi_{\alpha_i}}
               \\
                  uI \ar[r]^{\sigma}  
               & S
              }
$
\end{assumption}

In \cite{D} it is developed the case in which $\bb{S}$ is the category
of Sets, $\bb{I} = \{1\}$, and $u1 = 1 = \;$ the singleton set. 

\vspace{1ex}

\begin{proposition} [Under assumption \ref{classical}] \label{charstrictquasi}
Given a quasispace $(S,\, X)$, the following conditions are equivalent:

i) $(S,\, X)$ is a strict quasispace.

ii) $uI \mr{\sigma} S \; \in XI$ for all $I \in \bb{I}$ and $\sigma \in [uI, \, S]$.
\end{proposition}
\begin{proof}
Clearly  condition b) implies  $ii) \Rightarrow  i)$. The other implication follows by a straightforward application of conditions c) and a) on the family of all admissible maps. 
\end{proof}

Given any $S \in \bb{S}$, we can consider the quasispace structure generated by all maps 
$uI \to S$, $I \in \bb{I}$. This yields:

\begin{corollary}  [Under assumption \ref{classical}] \label{emptyfinalforqs}
Given any $S$ in $\bb{S}$, there exists a smallest strict quasispace structure on $S$, that we denote $S_{\bot_{\scriptstyle \ell}}$, and which is defined as follows:

\vspace{1ex}

$S_{\bot_{\scriptstyle \ell}} = (S,\, S_{\bot_{\scriptstyle \ell}}), \;\; where \;\;  uC \mr{\sigma} S  \, \in S_{\bot_{\scriptstyle \ell}} C\;\;$ if and only if $\;\;$ there exists 
$$
\xymatrix@R=5pt
               {
                \\
                C_i \to C \in \cc{T}_C, \;\;
                C_i \mr{}  I_i  \;\;  and \;\;  uI_i \mr{} S, \;\;\; such  \; that \;  
                }
\xymatrix@R=15pt@C=20pt
              {
                   uC_i \ar[r]  \ar[d]  
               &  uI_i  \ar[d]
               \\
                  uC \ar[r]^{\sigma}  
               & S
              }
$$
\qed\end{corollary} 

It follows then from \ref{topological} (\ref{chartop3}) and corollary
\ref{sandiareEfunctors} that the functor $s\bb{Q} \mr{q_s} \bb{S}$ is
topological. This can be seen directly as follows: 

Using proposition \ref{charstrictquasi} it is also immediate to check that given strict quasispaces  $(S_\alpha , \, X_\alpha)$ and any family $S \mr{\varphi_\alpha} S_\alpha$ in $\bb{S}$, the initial quasispace structure defined on $S$ in theorem \ref{Qistopological} by stipulating the equivalence in proposition \ref{charinitialQ} is strict. It follows:

\begin{corollary}  [Under assumption \ref{classical}] \label{initialinqs}
Given an initial family of quasispaces \mbox{$(S, \, X) \mr{\varphi_\alpha}  (S_\alpha ,\, X_\alpha)$,} if each  $(S_\alpha ,\, X_\alpha)$ is strict, then so it is  $(S, \, X)$.
\qed\end{corollary}

\begin{corollary}  [Under assumption \ref{classical}] \label{qsistopological}
The functor $s\bb{Q} \mr{q_s} \bb{S}$ is topological and  $s\bb{Q}$ is
closed under initial families in $\bb{Q}$.
\qed\end{corollary}

Warning: It does not follow that the functor $i$ is topological
(which is not).

\vspace{1ex}

From \ref{initialinqs} we have, in particular:

\begin{proposition}[Under assumption \ref{classical}]
Assumption \ref{a2} holds.
\qed \end{proposition}

Thus:

\begin{corollary}  [Under assumption \ref{classical}]
The functors 
$s\bb{Q} \mr{q_s} \bb{S}$ and $s\bb{Q} \mr{i} \bb{Q}$ are
$f$-regular, a category $s\bb{Q} \mr{q} \bb{S}$ of strict quasispaces
over a $f$-regular category  $\bb{S}$ is $f$-regular, and if $\bb{S}$
is a quasitopos, then so it is $s\bb{Q}$. Furthermore, the associate
quasispace functor situation described in \ref{aqfcontinued} holds.
\end{corollary}

\vspace{1ex}

The construction of the strict quasispace $S_{\bot_{\scriptstyle \ell}}$ in corollary \ref{emptyfinalforqs} can be generalized and provides a left adjoint for the inclusion $s\bb{Q} \mr{i} \bb{Q}$.

\begin{proposition} [Under assumption \ref{classical}]  \label{ihasaleftadjoint}  
The inclusion $s\bb{Q} \mr{i} \bb{Q}$ has a left adjoint  $\bb{Q} \mr{\ell} s\bb{Q}$,   
$\ell \dashv i$, such that $q_s \ell = q$ (notice that  $ (-)_{\bot_{\scriptstyle \ell}} = \ell (-)_\bot$).
\end{proposition}
\begin{proof}
We give an indication of the proof. The functor $\ell$ is easily understood: given a quasispace $(S,\, X)$, we abuse notation and denote $\ell(S, \, X) = (S, \, \ell X)$,  $X \subset \ell X$, where $\ell X$ is the quasispace generated by $X$ and all the arrows in $[uI, \, S]$, all $I \in \bb{I}$. Thus, 
$\ell(S, \, X) = (S, X \vee S_{\bot_{\scriptstyle \ell}})$.
\end{proof}

\end{document}